\documentclass[12pt,reqno]{amsart}

\usepackage{amsmath, epsfig, cite,color}
\usepackage{amssymb}
\usepackage{amsfonts}
\usepackage{latexsym}
\usepackage{graphicx}
\usepackage{algorithmic,algorithm}
\usepackage{amsthm,array,multirow,longtable}
\usepackage{tikz,caption}
\usetikzlibrary{decorations.pathreplacing}

\newtheorem{thm}{Theorem}[section]

\def\pf{\noindent {\it Proof.} }

\numberwithin{equation}{section}

\newcommand{\nc}{\newcommand}
\nc{\sha}{\mbox{\cyr X}}  
\font\cyr=wncyr10

\makeatletter \@addtoreset{equation}{section} \makeatother

\begin{document}
\allowdisplaybreaks
\title
{Combinatorics on lattice paths in strips}
\author[N.S.S. Gu]{Nancy S.S. Gu}
\address[N.S.S. Gu]{Center for Combinatorics, LPMC, Nankai
University, Tianjin 300071, P.R. China} \email{gu@nankai.edu.cn}
\author[H. Prodinger]{Helmut Prodinger}
\address[H. Prodinger]{Department of Mathematics, University of Stellenbosch,
7602 Stellenbosch, South Africa} \email{hproding@sun.ac.za}



\date{\today}

\begin{abstract}
For  lattice paths in strips which begin at $(0,0)$ and have only up steps $U: (i,j) \rightarrow (i+1,j+1)$ and down steps $D: (i,j)\rightarrow (i+1,j-1)$, let $A_{n,k}$ denote the set of  paths of length $n$ which start at $(0,0)$, end on heights $0$ or $-1$, and are contained
in the strip $-\lfloor\frac{k+1}{2}\rfloor \leq y \leq \lfloor\frac{k}{2}\rfloor$ of width $k$, and let $B_{n,k}$ denote the set of  paths of length $n$ which start at $(0,0)$ and are contained in the strip $0 \leq y \leq k$. We establish a bijection between $A_{n,k}$ and $B_{n,k}$.

The generating functions for the subsets of these two sets are discussed as well. Furthermore, we provide another bijection between $A_{n,3}$ and $B_{n,3}$ by translating the paths to two types of trees.
\end{abstract}

\maketitle

\section{Introduction}

Throughout this paper, we consider the lattice paths in strips which begin at $(0,0)$ and have only up steps $U: (i,j) \rightarrow (i+1,j+1)$ and down steps $D: (i,j)\rightarrow (i+1,j-1)$. Let $A_{n,k}$ denote the set of  paths of length $n$ which start at $(0,0)$, end on heights $0$ or $-1$, and are contained
in the strip $-\lfloor\frac{k+1}{2}\rfloor \leq y \leq \lfloor\frac{k}{2}\rfloor$ of width $k$.
Let $B_{n,k}$ denote the set of  paths of length $n$ which start at $(0,0)$ and are contained in the strip $0 \leq y \leq k$.
For example, we show $A_{4,3}$ and $B_{4,3}$ in Figures \ref{A-4-3} and \ref{B-4-3}, respectively.
\begin{figure}[h]
\begin{center}
\begin{tikzpicture}


\draw[->](0,-0.2)--(0,1.8);
\draw[->](-0.5,1)--(2.5,1);
\draw(0,1.5)--(2.5,1.5);
\draw(0,0.5)--(2.5,0.5);
\draw(0,0)--(2.5,0);
\draw(0,0)--(0.05,0);
\draw(0,0.5)--(0.05,0.5);
\draw(0,1)--(0.05,1);
\draw(0,1.5)--(0.05,1.5);
\node at (-0.3,0){\tiny$-2$};
\node at (-0.3,0.5){\tiny$-1$};
\node at (-0.2,0.85){\tiny$0$};
\node at (-0.2,1.5){\tiny$1$};
\node at (-0.2,1.8){\tiny$y$};
\draw(0.5,1)--(0.5,1.05);
\draw(1,1)--(1,1.05);
\draw(1.5,1)--(1.5,1.05);
\draw(2,1)--(2,1.05);
\node at (0.5,0.85){\tiny$1$};
\node at (1,0.85){\tiny$2$};
\node at (1.5,0.85){\tiny$3$};
\node at (2,0.85){\tiny$4$};
\node at (2.5,0.85){\tiny$x$};

\draw(0,1)--(0.5,0.5);
\draw(0.5,0.5)--(1,1);
\draw(1,1)--(1.5,0.5);
\draw(1.5,0.5)--(2,1);


\draw[->](3.5,-0.2)--(3.5,1.8);
\draw[->](3,1)--(6,1);
\draw(3.5,1.5)--(6,1.5);
\draw(3.5,0.5)--(6,0.5);
\draw(3.5,0)--(6,0);
\draw(3.5,0)--(3.55,0);
\draw(3.5,0.5)--(3.55,0.5);
\draw(3.5,1)--(3.55,1);
\draw(3.5,1.5)--(3.55,1.5);
\node at (3.2,0){\tiny$-2$};
\node at (3.2,0.5){\tiny$-1$};
\node at (3.3,0.85){\tiny$0$};
\node at (3.3,1.5){\tiny$1$};
\node at (3.3,1.8){\tiny$y$};
\draw(4,1)--(4,1.05);
\draw(4.5,1)--(4.5,1.05);
\draw(5,1)--(5,1.05);
\draw(5.5,1)--(5.5,1.05);
\node at (4,0.85){\tiny$1$};
\node at (4.5,0.85){\tiny$2$};
\node at (5,0.85){\tiny$3$};
\node at (5.5,0.85){\tiny$4$};
\node at (6,0.85){\tiny$x$};

\draw(3.5,1)--(4,0.5);
\draw(4,0.5)--(4.5,0);
\draw(4.5,0)--(5,0.5);
\draw(5,0.5)--(5.5,1);


\draw[->](0,2.3)--(0,4.3);
\draw[->](-0.5,3.5)--(2.5,3.5);
\draw(0,2.5)--(2.5,2.5);
\draw(0,3)--(2.5,3);
\draw(0,4)--(2.5,4);
\draw(0,2.5)--(0.05,2.5);
\draw(0,3)--(0.05,3);
\draw(0,3.5)--(0.05,3.5);
\draw(0,4)--(0.05,4);
\node at (-0.3,2.5){\tiny$-2$};
\node at (-0.3,3){\tiny$-1$};
\node at (-0.2,3.35){\tiny$0$};
\node at (-0.2,4){\tiny$1$};
\node at (-0.2,4.3){\tiny$y$};
\draw(0.5,3.5)--(0.5,3.55);
\draw(1,3.5)--(1,3.55);
\draw(1.5,3.5)--(1.5,3.55);
\draw(2,3.5)--(2,3.55);
\node at (0.5,3.35){\tiny$1$};
\node at (1,3.35){\tiny$2$};
\node at (1.5,3.35){\tiny$3$};
\node at (2,3.35){\tiny$4$};
\node at (2.5,3.35){\tiny$x$};

\draw(0,3.5)--(0.5,4);
\draw(0.5,4)--(1,3.5);
\draw(1,3.5)--(1.5,4);
\draw(1.5,4)--(2,3.5);


\draw[->](3.5,2.3)--(3.5,4.3);
\draw[->](3,3.5)--(6,3.5);
\draw(3.5,2.5)--(6,2.5);
\draw(3.5,3)--(6,3);
\draw(3.5,4)--(6,4);
\draw(3.5,2.5)--(3.55,2.5);
\draw(3.5,3)--(3.55,3);
\draw(3.5,3.5)--(3.55,3.5);
\draw(3.5,4)--(3.55,4);
\node at (3.2,2.5){\tiny$-2$};
\node at (3.2,3){\tiny$-1$};
\node at (3.3,3.35){\tiny$0$};
\node at (3.3,4){\tiny$1$};
\node at (3.3,4.3){\tiny$y$};
\draw(4,3.5)--(4,3.55);
\draw(4.5,3.5)--(4.5,3.55);
\draw(5,3.5)--(5,3.55);
\draw(5.5,3.5)--(5.5,3.55);
\node at (4,3.35){\tiny$1$};
\node at (4.5,3.35){\tiny$2$};
\node at (5,3.35){\tiny$3$};
\node at (5.5,3.35){\tiny$4$};
\node at (6,3.35){\tiny$x$};

\draw(3.5,3.5)--(4,4);
\draw(4,4)--(4.5,3.5);
\draw(4.5,3.5)--(5,3);
\draw(5,3)--(5.5,3.5);


\draw[->](7,2.3)--(7,4.3);
\draw[->](6.5,3.5)--(9.5,3.5);
\draw(7,2.5)--(9.5,2.5);
\draw(7,3)--(9.5,3);
\draw(7,4)--(9.5,4);
\draw(7,2.5)--(7.05,2.5);
\draw(7,3)--(7.05,3);
\draw(7,3.5)--(7.05,3.5);
\draw(7,4)--(7.05,4);
\node at (6.7,2.5){\tiny$-2$};
\node at (6.7,3){\tiny$-1$};
\node at (6.8,3.35){\tiny$0$};
\node at (6.8,4){\tiny$1$};
\node at (6.8,4.3){\tiny$y$};
\draw(7.5,3.5)--(7.5,3.55);
\draw(8,3.5)--(8,3.55);
\draw(8.5,3.5)--(8.5,3.55);
\draw(9,3.5)--(9,3.55);
\node at (7.5,3.35){\tiny$1$};
\node at (8,3.35){\tiny$2$};
\node at (8.5,3.35){\tiny$3$};
\node at (9,3.35){\tiny$4$};
\node at (9.5,3.35){\tiny$x$};

\draw(7,3.5)--(7.5,3);
\draw(7.5,3)--(8,3.5);
\draw(8,3.5)--(8.5,4);
\draw(8.5,4)--(9,3.5);

\end{tikzpicture}
\end{center}
\caption{The paths in $A_{4,3}$}\label{A-4-3}
\end{figure}
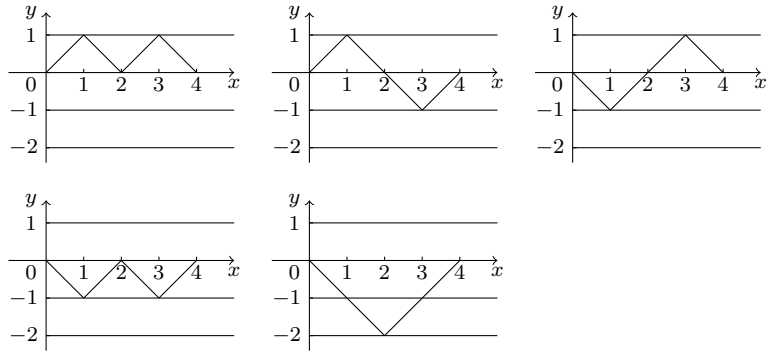

\begin{figure}[h]
\begin{center}
\begin{tikzpicture}


\draw[->](0,-0.3)--(0,1.8);
\draw[->](-0.5,0)--(2.5,0);
\draw(0,1.5)--(2.5,1.5);
\draw(0,1)--(2.5,1);
\draw(0,0.5)--(2.5,0.5);
\draw(0,0)--(0.05,0);
\draw(0,0.5)--(0.05,0.5);
\draw(0,1)--(0.05,1);
\draw(0,1.5)--(0.05,1.5);
\node at (-0.2,-0.15){\tiny$0$};
\node at (-0.2,0.5){\tiny$1$};
\node at (-0.2,1){\tiny$2$};
\node at (-0.2,1.5){\tiny$3$};
\node at (-0.2,1.8){\tiny$y$};
\draw(0.5,0)--(0.5,0.05);
\draw(1,0)--(1,0.05);
\draw(1.5,0)--(1.5,0.05);
\draw(2,0)--(2,0.05);
\node at (0.5,-0.15){\tiny$1$};
\node at (1,-0.15){\tiny$2$};
\node at (1.5,-0.15){\tiny$3$};
\node at (2,-0.15){\tiny$4$};
\node at (2.5,-0.15){\tiny$x$};

\draw(0,0)--(0.5,0.5);
\draw(0.5,0.5)--(1,0);
\draw(1,0)--(1.5,0.5);
\draw(1.5,0.5)--(2,1);


\draw[->](3.5,-0.3)--(3.5,1.8);
\draw[->](3,0)--(6,0);
\draw(3.5,1.5)--(6,1.5);
\draw(3.5,1)--(6,1);
\draw(3.5,0.5)--(6,0.5);
\draw(3.5,0)--(3.55,0);
\draw(3.5,0.5)--(3.55,0.5);
\draw(3.5,1)--(3.55,1);
\draw(3.5,1.5)--(3.55,1.5);
\node at (3.3,-0.15){\tiny$0$};
\node at (3.3,0.5){\tiny$1$};
\node at (3.3,1){\tiny$2$};
\node at (3.3,1.5){\tiny$3$};
\node at (3.3,1.8){\tiny$y$};
\draw(4,0)--(4,0.05);
\draw(4.5,0)--(4.5,0.05);
\draw(5,0)--(5,0.05);
\draw(5.5,0)--(5.5,0.05);
\node at (4,-0.15){\tiny$1$};
\node at (4.5,-0.15){\tiny$2$};
\node at (5,-0.15){\tiny$3$};
\node at (5.5,-0.15){\tiny$4$};
\node at (6,-0.15){\tiny$x$};

\draw(3.5,0)--(4,0.5);
\draw(4,0.5)--(4.5,0);
\draw(4.5,0)--(5,0.5);
\draw(5,0.5)--(5.5,0);


\draw[->](0,2.2)--(0,4.3);
\draw[->](-0.5,2.5)--(2.5,2.5);
\draw(0,4)--(2.5,4);
\draw(0,3.5)--(2.5,3.5);
\draw(0,3)--(2.5,3);
\draw(0,2.5)--(0.05,2.5);
\draw(0,3)--(0.05,3);
\draw(0,3.5)--(0.05,3.5);
\draw(0,4)--(0.05,4);
\node at (-0.2,2.35){\tiny$0$};
\node at (-0.2,3){\tiny$1$};
\node at (-0.2,3.5){\tiny$2$};
\node at (-0.2,4){\tiny$3$};
\node at (-0.2,4.3){\tiny$y$};
\draw(0.5,2.5)--(0.5,2.55);
\draw(1,2.5)--(1,2.55);
\draw(1.5,2.5)--(1.5,2.55);
\draw(2,2.5)--(2,2.55);
\node at (0.5,2.35){\tiny$1$};
\node at (1,2.35){\tiny$2$};
\node at (1.5,2.35){\tiny$3$};
\node at (2,2.35){\tiny$4$};
\node at (2.5,2.35){\tiny$x$};

\draw(0,2.5)--(0.5,3);
\draw(0.5,3)--(1,3.5);
\draw(1,3.5)--(1.5,4);
\draw(1.5,4)--(2,3.5);


\draw[->](3.5,2.2)--(3.5,4.3);
\draw[->](3,2.5)--(6,2.5);
\draw(3.5,4)--(6,4);
\draw(3.5,3.5)--(6,3.5);
\draw(3.5,3)--(6,3);
\draw(3.5,2.5)--(3.55,2.5);
\draw(3.5,3)--(3.55,3);
\draw(3.5,3.5)--(3.55,3.5);
\draw(3.5,4)--(3.55,4);
\node at (3.3,2.35){\tiny$0$};
\node at (3.3,3){\tiny$1$};
\node at (3.3,3.5){\tiny$2$};
\node at (3.3,4){\tiny$3$};
\node at (3.3,4.3){\tiny$y$};
\draw(4,2.5)--(4,2.55);
\draw(4.5,2.5)--(4.5,2.55);
\draw(5,2.5)--(5,2.55);
\draw(5.5,2.5)--(5.5,2.55);
\node at (4,2.35){\tiny$1$};
\node at (4.5,2.35){\tiny$2$};
\node at (5,2.35){\tiny$3$};
\node at (5.5,2.35){\tiny$4$};
\node at (6,2.35){\tiny$x$};

\draw(3.5,2.5)--(4,3);
\draw(4,3)--(4.5,3.5);
\draw(4.5,3.5)--(5,3);
\draw(5,3)--(5.5,3.5);


\draw[->](7,2.2)--(7,4.3);
\draw[->](6.5,2.5)--(9.5,2.5);
\draw(7,4)--(9.5,4);
\draw(7,3.5)--(9.5,3.5);
\draw(7,3)--(9.5,3);
\draw(7,2.5)--(7.05,2.5);
\draw(7,3)--(7.05,3);
\draw(7,3.5)--(7.05,3.5);
\draw(7,4)--(7.05,4);
\node at (6.8,2.35){\tiny$0$};
\node at (6.8,3){\tiny$1$};
\node at (6.8,3.5){\tiny$2$};
\node at (6.8,4){\tiny$3$};
\node at (6.8,4.3){\tiny$y$};
\draw(7.5,2.5)--(7.5,2.55);
\draw(8,2.5)--(8,2.55);
\draw(8.5,2.5)--(8.5,2.55);
\draw(9,2.5)--(9,2.55);
\node at (7.5,2.35){\tiny$1$};
\node at (8,2.35){\tiny$2$};
\node at (8.5,2.35){\tiny$3$};
\node at (9,2.35){\tiny$4$};
\node at (9.5,2.35){\tiny$x$};

\draw(7,2.5)--(7.5,3);
\draw(7.5,3)--(8,3.5);
\draw(8,3.5)--(8.5,3);
\draw(8.5,3)--(9,2.5);

\end{tikzpicture}
\end{center}
\caption{The paths in $B_{4,3}$}\label{B-4-3}
\end{figure}
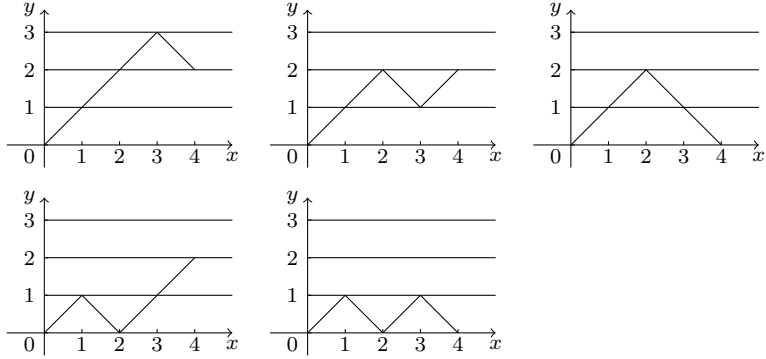

Using the viewpoint of the adjacency matrices of path graphs, Cigler \cite{cigler-1} presented that
\begin{equation*}
|A_{n,k}|=|B_{n,k}|=\sum_{j\in \mathbb{Z}}(-1)^j\binom{n}{\left\lfloor \frac{n+(k+2)j}{2} \right\rfloor}.
\end{equation*}
An interesting special case is that $A_{n,3}$ is enumerated by the Fibonacci number $F_{n+1}$ \cite[A000045]{OEIS}. For some other combinatorial objects which are related to the Fibonacci number, one can compare \cite{Barcucci,Czabarka,Deutsch-Prodinger,Florez}.  Cigler \cite{cigler-1} showed the following identity
\begin{equation*}
|A_{n,3}|=\sum_{j\in \mathbb{Z}}(-1)^j\binom{n}{\left\lfloor \frac{n+5j}{2} \right\rfloor}=F_{n+1}= \sum_{k=0}^{\lfloor\frac{n}{2}\rfloor}\binom{n-k}{k}
\end{equation*}
which is related to the first Rogers-Ramanujan identity
\begin{equation*}
\sum_{k=0}^{\infty}\frac{q^{k^2}}{(q;q)_k}=\frac{1}{(q;q^5)_{\infty}(q^4;q^5)_{\infty}},
\end{equation*}
where
$$(a;q)_k:=\prod_{i=0}^{k-1}(1-aq^{i})\qquad \text{and}\qquad (a;q)_{\infty}:=\prod_{i=0}^{\infty}(1-aq^{i}).$$

For another study of a similar type of  lattice paths, Panny and Prodinger \cite{Panny-Prodinger} considered the enumeration formulae for some types of expected height of lattice paths and their asymptotic equivalents.
Although Cigler \cite{cigler-1} derived the enumeration formulae for $A_{n,k}$ and $B_{n,k}$, he expected the existence of a simple bijection between $A_{n,k}$ and $B_{n,k}$. Prellberg \cite{prellberg} established a bijection between $A_{n,3}$ and $B_{n,3}$. Until now, Cigler's problem is still open for general $k$.
The main purpose of this paper is to solve the general case. The main theorem of this paper is stated as follows.
\begin{thm}\label{n-k}
There is a bijection between $A_{n,k}$ and $B_{n,k}$.
\end{thm}

This paper is organized as follows. In Section 2, we prove Theorem \ref{n-k}. In Section 3, we discuss the generating functions for the subsets of $A_{n,k}$ and $B_{n,k}$. In Section 4, we show another bijection between $A_{n,3}$ and $B_{n,3}$ by linking the paths to two types of trees.

\section{Main results}

To prove Theorem \ref{n-k}, we construct the bijection between $A_{2n,2k}$ (resp.\  $A_{2n, 2k+1}$, $A_{2n+1,2k}$, or $A_{2n+1,2k+1}$) and $B_{2n,2k}$ (resp.\  $B_{2n, 2k+1}$, $B_{2n+1,2k}$, or $B_{2n+1,2k+1}$), where
\begin{itemize}
\item $A_{2n,2k}$ (resp.\  $A_{2n+1,2k}$): the set of  paths of length $2n$ (resp.\  $2n+1$) which start at $(0,0)$, end on height $0$ (resp.\  $-1$), and are contained in the strip $-k \leq y \leq k$;
\item $A_{2n, 2k+1}$ (resp.\  $A_{2n+1,2k+1}$): the set of  paths of length $2n$ (resp.\  $2n+1$) which start at $(0,0)$, end on height $0$ (resp.\  $-1$), and are contained in the strip $-k-1 \leq y \leq k$;
\item $B_{2n,2k}$ (resp.\  $B_{2n+1,2k}$): the set of  paths of length $2n$ (resp.\  $2n+1$) which start at $(0,0)$ and are contained in the strip $0 \leq y \leq 2k$;
\item $B_{2n, 2k+1}$ (resp.\  $B_{2n+1,2k+1}$): the set of  paths of length $2n$ (resp.\  $2n+1$) which start at $(0,0)$ and are contained in the strip $0 \leq y \leq 2k+1$.
\end{itemize}

Notice that the paths in $B_{2n,2k}$ and $B_{2n, 2k+1}$ end on even heights, and the paths in $B_{2n+1,2k}$ and $B_{2n+1,2k+1}$ end on odd heights.

In the following, we prove four theorems. Then the proofs of these theorems imply Theorem \ref{n-k}. If there are no special requirements, we always look at the paths from left to right. For a given path $A$, let $\overline{A}$ denote the path which is obtained by flipping $A$ upside down. Conversely, if we flip the path $\overline{A}$ upside down, then we get the path $A$. For convenience, in all the pictures of this section, we only draw the lines parallel to the $x$-axis, and omit $x$-axis and $y$-axis.

\begin{thm}\label{2n-2k}
There is a bijection between $A_{2n,2k}$ and $B_{2n,2k}$.
\end{thm}
\pf We prove the theorem by induction on $k$.

When $k=1$, we separate $A_{2n,2}$ (resp.\  $B_{2n,2}$) into two subsets $A_{2n,2}^1$ and $A_{2n,2}^2$ (resp.\  $B_{2n,2}^1$ and $B_{2n,2}^2$).

$A_{2n,2}^1$ (resp.\  $A_{2n,2}^2$): the set of  paths in $A_{2n,2}$ which start with an up (resp.\  a down) step;

$B_{2n,2}^1$ (resp.\  $B_{2n,2}^2$): the set of  paths in $B_{2n,2}$ which end on height $2$ (resp.\  $0$).

For a path $P_1 \in A_{2n,2}^1$, let $a$ denote the first step of $P_1$, and let $A$ denote the sub-path after the step $a$  which starts with a down step and ends with an up or a down step. Here we take the picture in the upper-left corner of Figure 2.3 as an example to explain the meaning of the rectangle. The rectangle shows the scope of the strip in which the path $A$ exists. In this picture, it means that $A$ is in the strip $-1\leq y \leq 1$. The symbol
$\tikz [scale=0.2]\draw (0,1.2)to (1.2,0);$
in the upper-left corner of the rectangle represents the first down step of $A$ starting on height $1$, and the symbol
$\tikz [scale=0.2]\draw (0,1) to (1,0)to (0,-1);$ in the rectangle represents the last step of $A$ which is an up or a down step ending on height $0$. Since we do not care about the other steps of $A$, they are omitted in the rectangle.

We show the bijection between $A_{2n,2}^1$ (resp.\  $A_{2n,2}^2$) and $B_{2n,2}^1$ (resp.\  $B_{2n,2}^2$) in Figure \ref{2n-2}. The corresponding path $P_1' \in B_{2n,2}^1$ starts with $\overline{A}$ and ends with $a$, where $\overline{A}$ is obtained by flipping $A$ upside down. For a path $P_2 \in A_{2n,2}^2$, let $a$ denote the first step of $P_2$, and let $A$ denote the sub-path after the step $a$. Then we only change the order of $a$ and $A$ to obtain the corresponding path $P_2' \in B_{2n,2}^2$.

\begin{figure}[h]
\begin{center}
\begin{tikzpicture}


\draw(0,5)--(4,5);
\draw(0,4)--(4,4);
\draw(0,3)--(4,3);
\node at (-0.3,5){$1$};
\node at (-0.3,4){$0$};
\node at (-0.3,3){$-1$};

\draw(1,4)--(2,5);
\node at (1.3,4.6){$a$};

\draw(2,3)--(2,5);
\draw(3,3)--(3,5);
\draw(2,5)--(2.2,4.8);
\draw(3,4)--(2.8,4.2);
\draw(3,4)--(2.8,3.8);
\node at (2.5,4){$A$};

\node at (2,2.5){$A_{2n,2}^1$};

\node at (5,4){$\rightleftharpoons$};

\draw(6,5)--(10,5);
\draw(6,4)--(10,4);
\draw(6,3)--(10,3);
\node at (5.8,5){$2$};
\node at (5.8,4){$1$};
\node at (5.8,3){$0$};

\draw(8,4)--(9,5);
\node at (8.7,4.5){$a$};

\draw(7,3)--(7,5);
\draw(8,3)--(8,5);
\draw(7,3)--(7.2,3.2);
\draw(8,4)--(7.8,4.2);
\draw(8,4)--(7.8,3.8);
\node at (7.5,4){$\overline{A}$};

\node at (7.5,2.5){$B_{2n,2}^1$};


\draw(0,2)--(4,2);
\draw(0,1)--(4,1);
\draw(0,0)--(4,0);
\node at (-0.3,2){$1$};
\node at (-0.3,1){$0$};
\node at (-0.3,0){$-1$};

\draw(1,1)--(2,0);
\node at (1.3,0.4){$a$};

\draw(2,0)--(2,2);
\draw(3,0)--(3,2);
\draw(2,0)--(2.2,0.2);

\draw(3,1)--(2.8,1.2);
\draw(3,1)--(2.8,0.8);
\node at (2.5,1){$A$};

\node at (2,-0.5){$A_{2n,2}^2$};

\node at (5,1){$\rightleftharpoons$};

\draw(6,2)--(10,2);
\draw(6,1)--(10,1);
\draw(6,0)--(10,0);
\node at (5.8,2){$2$};
\node at (5.8,1){$1$};
\node at (5.8,0){$0$};

\draw(8,1)--(9,0);
\node at (8.7,0.7){$a$};

\draw(7,0)--(7,2);
\draw(8,0)--(8,2);
\draw(7,0)--(7.2,0.2);
\draw(8,1)--(7.8,1.2);
\draw(8,1)--(7.8,0.8);
\node at (7.5,1){$A$};

\node at (7.5,-0.5){$B_{2n,2}^2$};

\end{tikzpicture}
\end{center}
\caption{The bijection $\alpha_{2n,2}$ between $A_{2n,2}$ and $B_{2n,2}$}\label{2n-2}
\end{figure}

Suppose that there is a bijection $\alpha_{2n,2k-2}$ between $A_{2n,2k-2}$ and $B_{2n,2k-2}$ ($k\geq 2$). Let $A_{2n,2k}'$ denote the set of  paths in $A_{2n,2k}$ which touch the lines $y=k$ or $y=-k$.
Let $B_{2n,2k}'$ denote the set of  paths in $B_{2n,2k}$ such that in each path there exist some steps above the line  $y=2k-2$. Since $A_{2n,2k}$ (resp.\  $B_{2n,2k}$) can be divided into $A_{2n,2k-2}$ and $A_{2n,2k}'$ (resp.\  $B_{2n,2k-2}$ and $B_{2n,2k}'$), with the aid of the hypothesis, we only need to build the bijection between $A_{2n,2k}'$ and $B_{2n,2k}'$.

First, we separate $A_{2n,2k}'$ into three subsets $A_{2n,2k}^1$, $A_{2n,2k}^2$, and $A_{2n,2k}^3$, and separate $B_{2n,2k}'$ into $B_{2n,2k}^1$, $B_{2n,2k}^2$, and $B_{2n,2k}^3$.

$A_{2n,2k}^1$ (resp.\  $A_{2n,2k}^2$): the set of  paths in $A_{2n,2k}'$ which first arrive at height $k$, not height $-k$, such that in each path all the steps (resp.\  some steps) before the first point on height $k$ are above (resp.\  below) height $0$.

$A_{2n,2k}^3$: the set of  paths in $A_{2n,2k}'$ which first arrive at height $-k$, not height $k$.

$B_{2n,2k}^1$: the set of  paths in $B_{2n,2k}'$ which end on height $2k$.

$B_{2n,2k}^2$ (resp.\  $B_{2n,2k}^3$): the set of  paths in $B_{2n,2k}'$ which end on heights $0,2, \ldots$, or $2k-2$. Especially, in each path, starting from the last down step which ends on height $2k-2$, we go back from right to left until we arrive at height $k$ for the first time. This sub-path arrives (resp.\  doesn't arrive) at height $2k$.

We present the bijections $\alpha_{2n,2k}^i$ between $A_{2n,2k}^i$ and $B_{2n,2k}^i$ for $i=1,2,3$ in Figures \ref{2n-2k-1}--\ref{2n-2k-3}, respectively.

In Figure \ref{2n-2k-1}, for a given path in $A_{2n,2k}^1$, let $a$ denote the step which first arrives at height $k$. Then let $A$ (resp.\  $B$) denote the sub-path after (resp.\  before) the step $a$. To obtain the corresponding path in $B_{2n,2k}^1$, we first flip $A$ upside down to obtain $\overline{A}$. Then we put $\overline{A}$, $a$, and $B$ in order. To show the inverse map, for a given path in $B_{2n,2k}^1$, by finding the last up step $a$ starting on height $k$, we divide the path into three parts $\overline{A}$, $a$, and $B$. Then we flip $\overline{A}$ upside down to obtain $A$. Finally, we put $B$, $a$, and $A$ in order to obtain the corresponding path in $A_{2n,2k}^1$.

\begin{figure}[h]
\begin{center}
\begin{tikzpicture}


\draw(0,0)--(3.3,0);
\draw(0,1)--(3.3,1);
\draw(0,1.7)--(3.3,1.7);
\draw(0,2)--(3.3,2);

\node at (-0.5,0){\small$-k$};
\node at (-0.5,1){\small$0$};
\node at (-0.5,1.7){\small$k-1$};
\node at (-0.5,2){\small$k$};

\draw(0.5,1)--(0.5,1.7);
\draw(1.5,1)--(1.5,1.7);
\draw(0.5,1)--(0.7,1.2);
\draw(1.5,1.7)--(1.3,1.5);
\node at (1,1.3){$B$};

\draw(1.5,1.7)--(1.8,2);
\node at (1.5,1.85){$a$};

\draw(1.8,0)--(1.8,2);
\draw(2.8,0)--(2.8,2);
\draw(1.8,2)--(2,1.8);
\draw(2.8,1)--(2.6,1.2);
\draw(2.8,1)--(2.6,0.8);
\node at (2.3,1){$A$};

\node at (1.5,-0.5){$A_{2n,2k}^1$};

\node at (3.9,1){$\rightleftharpoons$};

\draw(5.5,2)--(8.8,2);
\draw(5.5,1.3)--(8.8,1.3);
\draw(5.5,1)--(8.8,1);
\draw(5.5,0)--(8.8,0);
\node at (5,0){\small$0$};
\node at (5,1){\small$k$};
\node at (5,1.3){\small$k+1$};
\node at (5,2){\small$2k$};

\draw(6,0)--(6,2);
\draw(7,0)--(7,2);
\draw(6,0)--(6.2,0.2);
\draw(7,1)--(6.8,1.2);
\draw(7,1)--(6.8,0.8);
\node at (6.5,1){$\overline{A}$};

\draw(7,1)--(7.3,1.3);
\node at (7.3,1.15){$a$};

\draw(7.3,1.3)--(7.3,2);
\draw(8.3,1.3)--(8.3,2);
\draw(7.3,1.3)--(7.5,1.5);
\draw(8.3,2)--(8.1,1.8);
\node at (7.8,1.7){$B$};

\node at (7,-0.5){$B_{2n,2k}^1$};

\end{tikzpicture}
\end{center}
\caption{The bijection $\alpha_{2n,2k}^1$ between $A_{2n,2k}^1$ and $B_{2n,2k}^1$}\label{2n-2k-1}
\end{figure}

In Figures \ref{2n-2k-2}, for a given path in $A_{2n,2k}^2$, let $a$ denote the step which first arrives at height $k$, and let $A$ denote the sub-path after the step $a$. Then let $b$ denote the last up step ending on height $0$ in the sub-path before the step $a$, and let $B$ denote the sub-path between $b$ and $a$. Next, from $b$, we go back from right to left to find the step $c$ when we first arrive at height $0$. Then the sub-path between $c$ and $b$ is denoted by $C$, and the sub-path before $c$ is denoted by $D$. So we divide the path into the parts $D$, $c$, $C$, $b$, $B$, $a$, and $A$. To show the map, first, since the sub-path $D$ is in the strip $-k+1\leq y \leq k-1$, we apply the hypothetical bijection on $D$, and flip the corresponding path upside down to obtain $D'$. Then by flipping $b$ and $A$ upside down, respectively, and putting $\overline{A}$, $a$, $B$, $\overline{b}$, $C$, $c$, and $D'$ in order, we obtain the corresponding path in $B_{2n,2k}^2$. Inversely, for a given path in $B_{2n,2k}^2$, let $c$ denote the last down step which ends on height $2k-2$. Then from $c$, we go back from right to left until we arrive at height $k$ for the first time. Let $a$ denote this step. For the sub-path between $a$ and $c$, from left to right, let $\overline{b}$ denote the last down step starting on height $2k$. Then we use $\overline{A}$, $B$, $C$, and $D'$ to denote the sub-paths which are separated by the steps $a$, $\overline{b}$, and $c$, respectively. To obtain the corresponding path in $A_{2n,2k}^2$, we first flip $D'$ upside down, and then apply the hypothetical bijection to get the path $D$. Then as shown in Figure \ref{2n-2k-2}, we complete the inverse map.

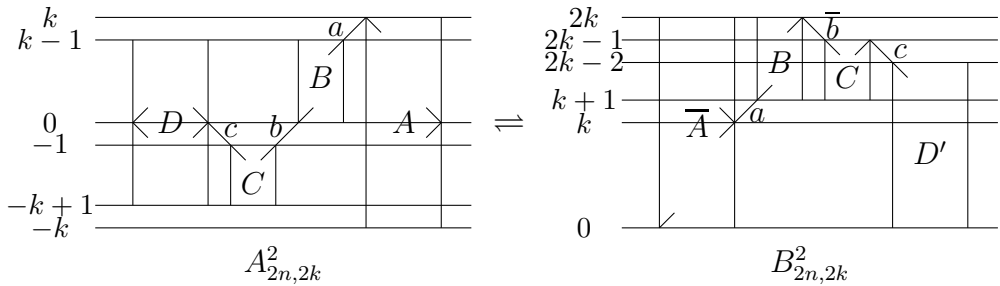
\begin{figure}[h]
\begin{center}
\begin{tikzpicture}


\draw(0,0)--(5,0);
\draw(0,0.3)--(5,0.3);
\draw(0,1.1)--(5,1.1);
\draw(0,1.4)--(5,1.4);
\draw(0,2.5)--(5,2.5);
\draw(0,2.8)--(5,2.8);
\node at (-0.6,0){\small$-k$};
\node at (-0.6,0.3){\small$-k+1$};
\node at (-0.6,1.1){\small$-1$};
\node at (-0.6,1.4){\small$0$};
\node at (-0.6,2.5){\small$k-1$};
\node at (-0.6,2.8){\small$k$};

\draw(0.5,0.3)--(0.5,2.5);
\draw(1.5,0.3)--(1.5,2.5);
\draw(0.5,1.4)--(0.7,1.2);
\draw(0.5,1.4)--(0.7,1.6);
\draw(1.5,1.4)--(1.3,1.2);
\draw(1.5,1.4)--(1.3,1.6);
\node at (1,1.4){$D$};

\draw(1.5,1.4)--(1.8,1.1);
\node at (1.8,1.3){$c$};

\draw(1.8,1.1)--(1.8,0.3);
\draw(2.4,1.1)--(2.4,0.3);
\draw(1.8,1.1)--(2,0.9);
\draw(2.4,1.1)--(2.2,0.9);
\node at (2.1,0.6){$C$};

\draw(2.4,1.1)--(2.7,1.4);
\node at (2.4,1.35){$b$};

\draw(2.7,1.4)--(2.7,2.5);
\draw(3.3,1.4)--(3.3,2.5);
\draw(2.7,1.4)--(2.9,1.6);
\draw(3.3,2.5)--(3.1,2.3);
\node at (3,2){$B$};

\draw(3.3,2.5)--(3.6,2.8);
\node at (3.2,2.65){$a$};

\draw(3.6,2.8)--(3.6,0);
\draw(4.6,2.8)--(4.6,0);
\draw(3.6,2.8)--(3.8,2.6);
\draw(4.6,1.4)--(4.4,1.2);
\draw(4.6,1.4)--(4.4,1.6);
\node at (4.1,1.4){$A$};

\node at (2.5,-0.5){$A_{2n,2k}^2$};

\node at (5.5,1.4){$\rightleftharpoons$};

\draw(7,0)--(12,0);
\draw(7,1.4)--(12,1.4);
\draw(7,1.7)--(12,1.7);
\draw(7,2.2)--(12,2.2);
\draw(7,2.5)--(12,2.5);
\draw(7,2.8)--(12,2.8);
\node at (6.5,0){\small$0$};
\node at (6.5,1.4){\small$k$};
\node at (6.5,1.7){\small$k+1$};
\node at (6.5,2.2){\small$2k-2$};
\node at (6.5,2.5){\small$2k-1$};
\node at (6.5,2.8){\small$2k$};

\draw(7.5,0)--(7.5,2.8);
\draw(8.5,0)--(8.5,2.8);
\draw(7.5,0)--(7.7,0.2);
\draw(8.5,1.4)--(8.3,1.2);
\draw(8.5,1.4)--(8.3,1.6);
\node at (8,1.4){$\overline{A}$};

\draw(8.5,1.4)--(8.8,1.7);
\node at (8.8,1.5){$a$};

\draw(8.8,1.7)--(8.8,2.8);
\draw(9.4,1.7)--(9.4,2.8);
\draw(8.8,1.7)--(9,1.9);
\draw(9.4,2.8)--(9.2,2.6);
\node at (9.1,2.2){$B$};

\draw(9.4,2.8)--(9.7,2.5);
\node at (9.8,2.7){$\overline{b}$};

\draw(9.7,2.5)--(9.7,1.7);
\draw(10.3,2.5)--(10.3,1.7);
\draw(9.7,2.5)--(9.9,2.3);
\draw(10.3,2.5)--(10.1,2.3);
\node at (10,2){$C$};

\draw(10.3,2.5)--(10.6,2.2);
\node at (10.7,2.35){$c$};

\draw(10.6,2.2)--(10.6,0);
\draw(11.6,2.2)--(11.6,0);
\draw(10.6,2.2)--(10.8,2);
\node at (11.1,1){$D'$};

\node at (9.5,-0.5){$B_{2n,2k}^2$};

\end{tikzpicture}
\end{center}
\caption{The bijection $\alpha_{2n,2k}^2$ between $A_{2n,2k}^2$ and $B_{2n,2k}^2$}\label{2n-2k-2}
\end{figure}

In Figures \ref{2n-2k-3}, for a given path in $A_{2n,2k}^3$, let $a$ denote the first down step ending on height $-k$, and let $A$ denote the sub-path after the step $a$. Then from $a$, we go back from right to left until we find the first step ending on height $0$, and let $b$ denote this step. Alternatively, we use $B$ to denote the sub-path between $b$ and $a$, and let $C$ denote the sub-path before $b$. To show the map from $A_{2n,2k}^3$ to $B_{2n,2k}^3$, first, we apply the hypothetical bijection on $C$, and then flip the corresponding path upside down to obtain $C'$. Next, by flipping $a$ and $B$ upside down, respectively, and then putting $A$, $\overline{a}$, $\overline{B}$, $b$, and $C'$ in order, we get the corresponding path in $B_{2n,2k}^3$. Inversely, for a given path in $B_{2n,2k}^3$, first, let $b$ denote the last down step ending on height $2k-2$, and let $C'$ denote the sub-path after the step $b$. Then from $b$, we go back from right to left until we first arrive at height $k$. Let $\overline{a}$ denote this step. Furthermore, let $\overline{B}$ denote the sub-path between $\overline{a}$ and $b$, and let $A$ denote the sub-path before $\overline{a}$. To obtain the corresponding path in $A_{2n,2k}^3$, we first flip $C'$ upside down, and then apply the hypothetical bijection to obtain $C$. Next, by flipping $\overline{a}$ and $\overline{B}$ upside down, respectively, and then putting $C$, $b$, $B$, $a$, and $A$ in order, we complete the inverse map.

\begin{figure}[h]
\begin{center}
\begin{tikzpicture}


\draw(0,0)--(4.6,0);
\draw(0,0.3)--(4.6,0.3);
\draw(0,1.1)--(4.6,1.1);
\draw(0,1.4)--(4.6,1.4);
\draw(0,2.5)--(4.6,2.5);
\draw(0,2.8)--(4.6,2.8);
\node at (-0.6,0){\small$-k$};
\node at (-0.6,0.3){\small$-k+1$};
\node at (-0.6,1.1){\small$-1$};
\node at (-0.6,1.4){\small$0$};
\node at (-0.6,2.5){\small$k-1$};
\node at (-0.6,2.8){\small$k$};

\draw(0.5,0.3)--(0.5,2.5);
\draw(1.5,0.3)--(1.5,2.5);
\draw(0.5,1.4)--(0.7,1.2);
\draw(0.5,1.4)--(0.7,1.6);
\draw(1.5,1.4)--(1.3,1.2);
\draw(1.5,1.4)--(1.3,1.6);
\node at (1,1.4){$C$};

\draw(1.5,1.4)--(1.8,1.1);
\node at (1.8,1.3){$b$};

\draw(1.8,1.1)--(1.8,0.3);
\draw(2.8,1.1)--(2.8,0.3);
\draw(1.8,1.1)--(2,0.9);
\draw(2.8,0.3)--(2.6,0.5);
\node at (2.3,0.6){$B$};

\draw(2.8,0.3)--(3.1,0);
\node at (2.8,0.1){$a$};

\draw(3.1,2.8)--(3.1,0);
\draw(4.1,2.8)--(4.1,0);
\draw(3.1,0)--(3.3,0.2);
\draw(4.1,1.4)--(3.9,1.2);
\draw(4.1,1.4)--(3.9,1.6);
\node at (3.6,1.4){$A$};

\node at (2.2,-0.5){$A_{2n,2k}^3$};

\node at (5.5,1.4){$\rightleftharpoons$};

\draw(7,0)--(11.6,0);
\draw(7,1.4)--(11.6,1.4);
\draw(7,1.7)--(11.6,1.7);
\draw(7,2.2)--(11.6,2.2);
\draw(7,2.5)--(11.6,2.5);
\draw(7,2.8)--(11.6,2.8);
\node at (6.5,0){\small$0$};
\node at (6.5,1.4){\small$k$};
\node at (6.5,1.7){\small$k+1$};
\node at (6.5,2.2){\small$2k-2$};
\node at (6.5,2.5){\small$2k-1$};
\node at (6.5,2.8){\small$2k$};

\draw(7.5,0)--(7.5,2.8);
\draw(8.5,0)--(8.5,2.8);
\draw(7.5,0)--(7.7,0.2);
\draw(8.5,1.4)--(8.3,1.2);
\draw(8.5,1.4)--(8.3,1.6);
\node at (8,1.4){$A$};

\draw(8.5,1.4)--(8.8,1.7);
\node at (8.8,1.5){$\overline{a}$};

\draw(8.8,1.7)--(8.8,2.5);
\draw(9.8,1.7)--(9.8,2.5);
\draw(8.8,1.7)--(9,1.9);
\draw(9.8,2.5)--(9.6,2.3);
\node at (9.3,2.15){$\overline{B}$};

\draw(9.8,2.5)--(10.1,2.2);
\node at (10.2,2.45){$b$};

\draw(10.1,2.2)--(10.1,0);
\draw(11.1,2.2)--(11.1,0);
\draw(10.1,2.2)--(10.3,2);
\node at (10.6,1){$C'$};

\node at (9.5,-0.5){$B_{2n,2k}^3$};

\end{tikzpicture}
\end{center}
\caption{The bijection $\alpha_{2n,2k}^3$ between $A_{2n,2k}^3$ and $B_{2n,2k}^3$}\label{2n-2k-3}
\end{figure}
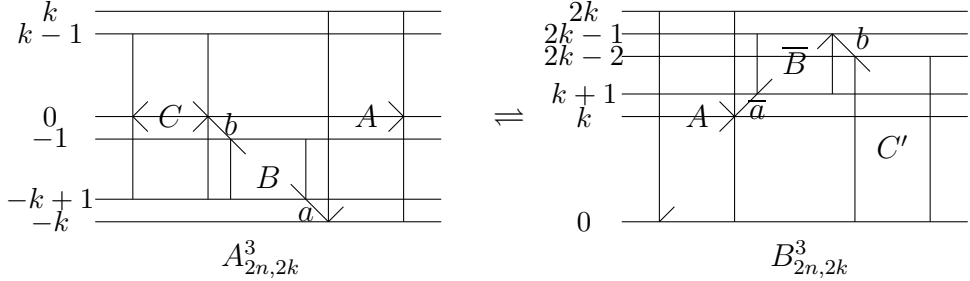

Notice that when $k=2$, the sub-path $B$ doesn't exist in Figure \ref{2n-2k-3}. So for any path in $B_{2n,4}^3$ the down step $b$ starting on height $3$  must follow $\overline{a}$.

In conclusion, combining the bijections $\alpha_{2n,2k}^i$ for $i=1,2,3$ yields the bijection $\alpha_{2n,2k}$ between $A_{2n,2k}$ and $B_{2n,2k}$.
\qed

We use the following example to explain the bijection in Theorem \ref{2n-2k}. First, we give a path $P\in A_{26,6}$ on the left in Figure \ref{ex-1}. Since $P$ first arrives at height $-3$, not height $3$, we have $P \in A_{26,6}^3$. Let $a$ denote the first down step ending on height $-3$. From $a$, we go back from right to left until we first arrive at height $0$, and let $b$ denote this step. Then with the aid of the steps $a$ and $b$, we can mark the other three sub-paths by $A$, $B$, and $C$, respectively. Applying the bijection $\alpha_{26,6}^3$ on $P$, we derive the corresponding path on the right in Figure \ref{ex-1}. Since $C$ is a path in the strip $-2\leq y\leq 2$, we need to use the bijection $\alpha_{14,4}$, and then flip the corresponding path upside down to derive the path $C'$. Here we only use the rectangle with a down step to represent $C'$, and explain the detailed map from $C$ to $C'$ in the following.

\begin{figure}[h]
\begin{center}
\begin{tikzpicture}

\draw (0,0)--(8.4,0);
\draw (0,0.3)--(8.4,0.3);
\draw (0,0.6)--(8.4,0.6);
\draw (0,0.9)--(8.4,0.9);
\draw (0,1.2)--(8.4,1.2);
\draw (0,1.5)--(8.4,1.5);
\draw (0,1.8)--(8.4,1.8);
\node at (-0.3,0){\tiny$-3$};
\node at (-0.3,0.3){\tiny$-2$};
\node at (-0.3,0.6){\tiny$-1$};
\node at (-0.2,0.9){\tiny$0$};
\node at (-0.2,1.2){\tiny$1$};
\node at (-0.2,1.5){\tiny$2$};
\node at (-0.2,1.8){\tiny$3$};

\draw (0.3,0.9)--(0.6,0.6);
\draw (0.6,0.6)--(0.9,0.9);
\draw (0.9,0.9)--(1.2,1.2);
\draw (1.2,1.2)--(1.5,0.9);
\draw (1.5,0.9)--(1.8,0.6);
\draw (1.8,0.6)--(2.1,0.9);
\draw (2.1,0.9)--(2.4,1.2);
\draw (2.4,1.2)--(2.7,1.5);
\draw (2.7,1.5)--(3,1.2);
\draw (3,1.2)--(3.3,1.5);
\draw (3.3,1.5)--(3.6,1.2);
\draw (3.6,1.2)--(3.9,0.9);
\draw (3.9,0.9)--(4.2,0.6);
\draw (4.2,0.6)--(4.5,0.9);
\draw (4.5,0.9)--(4.8,0.6);
\draw (4.8,0.6)--(5.1,0.3);
\draw (5.1,0.3)--(5.4,0.6);
\draw (5.4,0.6)--(5.7,0.3);
\draw (5.7,0.3)--(6,0);
\draw (6,0)--(6.3,0.3);
\draw (6.3,0.3)--(6.6,0);
\draw (6.6,0)--(6.9,0.3);
\draw (6.9,0.3)--(7.2,0.6);
\draw (7.2,0.6)--(7.5,0.9);
\draw (7.5,0.9)--(7.8,1.2);
\draw (7.8,1.2)--(8.1,0.9);

\draw (0.3,1.5)--(0.3,0.3);
\draw (4.5,1.5)--(4.5,0.3);

\draw (4.8,0.6)--(4.8,0.3);
\draw (5.7,0.6)--(5.7,0.3);

\draw (6,1.8)--(6,0);
\draw (8.1,1.8)--(8.1,0);

\node at (3,1.05){\tiny$C$};
\node at (4.75,0.8){\tiny$b$};
\node at (5.4,0.42){\tiny$B$};
\node at (5.7,0.15){\tiny$a$};
\node at (7,0.75){\tiny$A$};

\node at (8.8,0.9){$\rightleftharpoons$};
\node at (8.9,1.25){\tiny$\alpha_{26,6}^3$};

\draw (9.6,0)--(14.8,0);
\draw (9.6,0.3)--(14.8,0.3);
\draw (9.6,0.6)--(14.8,0.6);
\draw (9.6,0.9)--(14.8,0.9);
\draw (9.6,1.2)--(14.8,1.2);
\draw (9.6,1.5)--(14.8,1.5);
\draw (9.6,1.8)--(14.8,1.8);
\node at (9.4,0){\tiny$0$};
\node at (9.4,0.3){\tiny$1$};
\node at (9.4,0.6){\tiny$2$};
\node at (9.4,0.9){\tiny$3$};
\node at (9.4,1.2){\tiny$4$};
\node at (9.4,1.5){\tiny$5$};
\node at (9.4,1.8){\tiny$6$};

\draw (9.9,0)--(10.2,0.3);
\draw (10.2,0.3)--(10.5,0);
\draw (10.5,0)--(10.8,0.3);
\draw (10.8,0.3)--(11.1,0.6);
\draw (11.1,0.6)--(11.4,0.9);
\draw (11.4,0.9)--(11.7,1.2);
\draw (11.7,1.2)--(12,0.9);
\draw (12,0.9)--(12.3,1.2);
\draw (12.3,1.2)--(12.6,1.5);
\draw (12.6,1.5)--(12.9,1.2);
\draw (12.9,1.2)--(13.2,1.5);
\draw (13.2,1.5)--(13.5,1.2);
\draw (13.5,1.2)--(13.7,1);
\draw (13.5,1.2)--(13.5,0);
\draw (14.5,1.2)--(14.5,0);

\draw (9.9,1.8)--(9.9,0);
\draw (12,1.8)--(12,0);

\draw (12.3,1.5)--(12.3,1.2);
\draw (13.2,1.5)--(13.2,1.2);

\node at (10.9,0.75){\tiny$A$};
\node at (12.3,1.02){\tiny$\overline{a}$};
\node at (12.9,1.44){\tiny$\overline{B}$};
\node at (13.5,1.38){\tiny$b$};
\node at (14,0.75){\tiny$C'$};

\end{tikzpicture}
\end{center}
\caption{An example for the bijection $\alpha_{2n,2k}$}\label{ex-1}
\end{figure}

In Figure \ref{ex-2}, we show how to obtain $\overline{C'}$ from $C$. For the path $C \in A_{14,4}^2$, based on the illustration of $\alpha_{2n,2k}^2$ in Theorem \ref{2n-2k}, we mark the sub-paths by $D$, $c$, $b$, $B$, $a$ and $A$ in the path on the left in Figure \ref{ex-2}, respectively. To obtain the corresponding path $\overline{C'}$, we apply the bijection $\alpha_{14,4}^2$ to get the path on the right in Figure \ref{ex-2}. Especially, we apply the bijection $\alpha_{4,2}$ on $D$, and then flip the corresponding path upside down to get $D'$.

\begin{figure}[h]
\begin{center}
\begin{tikzpicture}

\draw (0,0)--(4.8,0);
\draw (0,0.3)--(4.8,0.3);
\draw (0,0.6)--(4.8,0.6);
\draw (0,0.9)--(4.8,0.9);
\draw (0,1.2)--(4.8,1.2);

\node at (-0.3,0){\tiny$-2$};
\node at (-0.3,0.3){\tiny$-1$};
\node at (-0.2,0.6){\tiny$0$};
\node at (-0.2,0.9){\tiny$1$};
\node at (-0.2,1.2){\tiny$2$};

\draw (0.3,0.6)--(0.6,0.3);
\draw (0.6,0.3)--(0.9,0.6);
\draw (0.9,0.6)--(1.2,0.9);
\draw (1.2,0.9)--(1.5,0.6);
\draw (1.5,0.6)--(1.8,0.3);
\draw (1.8,0.3)--(2.1,0.6);
\draw (2.1,0.6)--(2.4,0.9);
\draw (2.4,0.9)--(2.7,1.2);
\draw (2.7,1.2)--(3,0.9);
\draw (3,0.9)--(3.3,1.2);
\draw (3.3,1.2)--(3.6,0.9);
\draw (3.6,0.9)--(3.9,0.6);
\draw (3.9,0.6)--(4.2,0.3);
\draw (4.2,0.3)--(4.5,0.6);

\draw (0.3,0.9)--(0.3,0.3);
\draw (1.5,0.9)--(1.5,0.3);

\draw (2.1,0.9)--(2.1,0.6);
\draw (2.4,0.9)--(2.4,0.6);

\draw (2.7,1.2)--(2.7,0);
\draw (4.5,1.2)--(4.5,0);

\node at (1.1,0.45){\tiny$D$};
\node at (1.59,0.39){\tiny$c$};
\node at (2.05,0.41){\tiny$b$};
\node at (2.25,0.75){\tiny$B$};
\node at (2.44,1.08){\tiny$a$};
\node at (3.3,0.75){\tiny$A$};

\node at (2.4,-0.5){\tiny$C$};

\node at (5.8,0.95){\tiny$\alpha_{14,4}^2$};
\node at (5.75,0.6){$\rightleftharpoons$};

\draw (7,0)--(11.8,0);
\draw (7,0.3)--(11.8,0.3);
\draw (7,0.6)--(11.8,0.6);
\draw (7,0.9)--(11.8,0.9);
\draw (7,1.2)--(11.8,1.2);

\node at (6.8,0){\tiny$0$};
\node at (6.8,0.3){\tiny$1$};
\node at (6.8,0.6){\tiny$2$};
\node at (6.8,0.9){\tiny$3$};
\node at (6.8,1.2){\tiny$4$};

\draw (7.3,0)--(7.6,0.3);
\draw (7.6,0.3)--(7.9,0);
\draw (7.9,0)--(8.2,0.3);
\draw (8.2,0.3)--(8.5,0.6);
\draw (8.5,0.6)--(8.8,0.9);
\draw (8.8,0.9)--(9.1,0.6);
\draw (9.1,0.6)--(9.4,0.9);
\draw (9.4,0.9)--(9.7,1.2);
\draw (9.7,1.2)--(10,0.9);
\draw (10,0.9)--(10.3,0.6);
\draw (10.3,0.6)--(10.6,0.3);
\draw (10.6,0.3)--(10.9,0);
\draw (10.9,0)--(11.2,0.3);
\draw (11.2,0.3)--(11.5,0.6);

\draw (7.3,1.2)--(7.3,0);
\draw (9.1,1.2)--(9.1,0);

\draw (9.4,1.2)--(9.4,0.9);
\draw (9.7,1.2)--(9.7,0.9);

\draw (10.3,0.6)--(10.3,0);
\draw (11.5,0.6)--(11.5,0);

\node at (7.9,0.45){\tiny$\overline{A}$};
\node at (9.35,0.73){\tiny$a$};
\node at (9.55,1.05){\tiny$B$};
\node at (10.05,1.05){\tiny$\overline{b}$};
\node at (10.3,0.75){\tiny$c$};
\node at (10.9,0.45){\tiny$D'$};

\node at (9.4,-0.5){\tiny$\overline{C'}$};

\end{tikzpicture}
\end{center}
\caption{The bijection between $C$ and $\overline{C'}$}\label{ex-2}
\end{figure}

\noindent Combining Figures \ref{ex-1} and \ref{ex-2} yields the final path $P' \in B_{26,6}$ in Figure \ref{ex-3}.

\begin{figure}[h]
\begin{center}
\begin{tikzpicture}

\draw (0,0)--(8.4,0);
\draw (0,0.3)--(8.4,0.3);
\draw (0,0.6)--(8.4,0.6);
\draw (0,0.9)--(8.4,0.9);
\draw (0,1.2)--(8.4,1.2);
\draw (0,1.5)--(8.4,1.5);
\draw (0,1.8)--(8.4,1.8);

\node at (-0.2,0){\tiny$0$};
\node at (-0.2,0.3){\tiny$1$};
\node at (-0.2,0.6){\tiny$2$};
\node at (-0.2,0.9){\tiny$3$};
\node at (-0.2,1.2){\tiny$4$};
\node at (-0.2,1.5){\tiny$5$};
\node at (-0.2,1.8){\tiny$6$};

\draw (0.3,0)--(0.6,0.3);
\draw (0.6,0.3)--(0.9,0);
\draw (0.9,0)--(1.2,0.3);
\draw (1.2,0.3)--(1.5,0.6);
\draw (1.5,0.6)--(1.8,0.9);
\draw (1.8,0.9)--(2.1,1.2);
\draw (2.1,1.2)--(2.4,0.9);
\draw (2.4,0.9)--(2.7,1.2);
\draw (2.7,1.2)--(3,1.5);
\draw (3,1.5)--(3.3,1.2);
\draw (3.3,1.2)--(3.6,1.5);
\draw (3.6,1.5)--(3.9,1.2);
\draw (3.9,1.2)--(4.2,0.9);
\draw (4.2,0.9)--(4.5,1.2);
\draw (4.5,1.2)--(4.8,0.9);
\draw (4.8,0.9)--(5.1,0.6);
\draw (5.1,0.6)--(5.4,0.3);
\draw (5.4,0.3)--(5.7,0.6);
\draw (5.7,0.6)--(6,0.3);
\draw (6,0.3)--(6.3,0);
\draw (6.3,0)--(6.6,0.3);
\draw (6.6,0.3)--(6.9,0.6);
\draw (6.9,0.6)--(7.2,0.9);
\draw (7.2,0.9)--(7.5,1.2);
\draw (7.5,1.2)--(7.8,0.9);
\draw (7.8,0.9)--(8.1,0.6);

\node at (2.7,1.05){\tiny$\overline{a}$};
\node at (3.9,1.37){\tiny$b$};

\end{tikzpicture}
\end{center}
\caption{The path $P'\in B_{26,6}$}\label{ex-3}
\end{figure}

Inversely, for the given path $P'\in B_{26,6}$ in Figure \ref{ex-3}, we first find the last down step $b$ ending on height $4$. Then from $b$, we go back from right to left to find the step $\overline{a}$ when we first arrive at height $3$. Since the sub-path between $\overline{a}$ and $b$ doesn't arrive at height $6$, we have $P'\in B_{26,6}^3$. Then based on the bijections in Figures \ref{ex-1} and \ref{ex-2}, we can derive the path $P \in A_{26,6}$.

\begin{thm}
There is a bijection between $A_{2n,2k+1}$ and $B_{2n,2k+1}$.
\end{thm}
\pf  The set $A_{2n,2k+1}$ (resp.\  $B_{2n,2k+1}$) can be divided into two subsets $A_{2n,2k+1}^1$ and $A_{2n,2k+1}^2$
(resp.\  $B_{2n,2k+1}^1$ and $B_{2n,2k+1}^2$).

$A_{2n,2k+1}^1$: the set of  paths in $A_{2n,2k+1}$ which are in the strip $-k \leq y \leq k$.

$A_{2n,2k+1}^2$: the set of  paths in $A_{2n,2k+1}$ which arrive at height $-k-1$.

$B_{2n,2k+1}^1$: the set of  paths in $B_{2n,2k+1}$ which are in the strip $0 \leq y \leq 2k$.

$B_{2n,2k+1}^2$: the set of  paths in $B_{2n,2k+1}$ which arrive at height $2k+1$.

\begin{figure}[h]
\begin{center}
\begin{tikzpicture}


\draw(0,0)--(4.6,0);
\draw(0,0.3)--(4.6,0.3);
\draw(0,1.4)--(4.6,1.4);
\draw(0,1.7)--(4.6,1.7);
\draw(0,3.1)--(4.6,3.1);
\node at (-0.6,0){\small$-k-1$};
\node at (-0.6,0.3){\small$-k$};
\node at (-0.6,1.4){\small$-1$};
\node at (-0.6,1.7){\small$0$};
\node at (-0.6,3.1){\small$k$};

\draw(0.5,0.3)--(0.5,3.1);
\draw(1.5,0.3)--(1.5,3.1);
\draw(0.5,1.7)--(0.7,1.5);
\draw(0.5,1.7)--(0.7,1.9);
\draw(1.5,1.7)--(1.3,1.5);
\draw(1.5,1.7)--(1.3,1.9);
\node at (1,1.7){$C$};

\draw(1.5,1.7)--(1.8,1.4);
\node at (1.8,1.65){$b$};

\draw(1.8,1.4)--(1.8,0.3);
\draw(2.8,1.4)--(2.8,0.3);
\draw(1.8,1.4)--(2,1.2);
\draw(2.8,0.3)--(2.6,0.5);
\node at (2.3,0.8){$B$};

\draw(2.8,0.3)--(3.1,0);
\node at (2.8,0.1){$a$};

\draw(3.1,3.1)--(3.1,0);
\draw(4.1,3.1)--(4.1,0);
\draw(3.1,0)--(3.3,0.2);
\draw(4.1,1.7)--(3.9,1.5);
\draw(4.1,1.7)--(3.9,1.9);
\node at (3.6,1.4){$A$};

\node at (2.2,-0.5){$A_{2n,2k+1}^2$};

\node at (5.5,1.4){$\rightleftharpoons$};

\draw(7,0)--(11.6,0);
\draw(7,1.7)--(11.6,1.7);
\draw(7,2)--(11.6,2);
\draw(7,2.5)--(11.6,2.5);
\draw(7,2.8)--(11.6,2.8);
\draw(7,3.1)--(11.6,3.1);
\node at (6.5,0){\small$0$};
\node at (6.5,1.7){\small$k+1$};
\node at (6.5,2){\small$k+2$};
\node at (6.5,2.5){\small$2k-1$};
\node at (6.5,2.8){\small$2k$};
\node at (6.5,3.1){\small$2k+1$};

\draw(7.5,0)--(7.5,3.1);
\draw(8.5,0)--(8.5,3.1);
\draw(7.5,0)--(7.7,0.2);
\draw(8.5,1.7)--(8.3,1.5);
\draw(8.5,1.7)--(8.3,1.9);
\node at (8,1.4){$A$};

\draw(8.5,1.7)--(8.8,2);
\node at (8.8,1.8){$\overline{a}$};

\draw(8.8,2)--(8.8,3.1);
\draw(9.8,2)--(9.8,3.1);
\draw(8.8,2)--(9,2.2);
\draw(9.8,3.1)--(9.6,2.9);
\node at (9.3,2.5){$\overline{B}$};

\draw(9.8,3.1)--(10.1,2.8);
\node at (10.2,3){$b$};

\draw(10.1,2.8)--(10.1,0);
\draw(11.1,2.8)--(11.1,0);
\draw(10.1,2.8)--(10.3,2.6);
\node at (10.6,1.4){$C'$};

\node at (9.5,-0.5){$B_{2n,2k+1}^2$};

\end{tikzpicture}
\end{center}
\caption{The bijection $\alpha_{2n,2k+1}^2$ between $A_{2n,2k+1}^2$ and $B_{2n,2k+1}^2$}\label{2n-2k+1-2}
\end{figure}
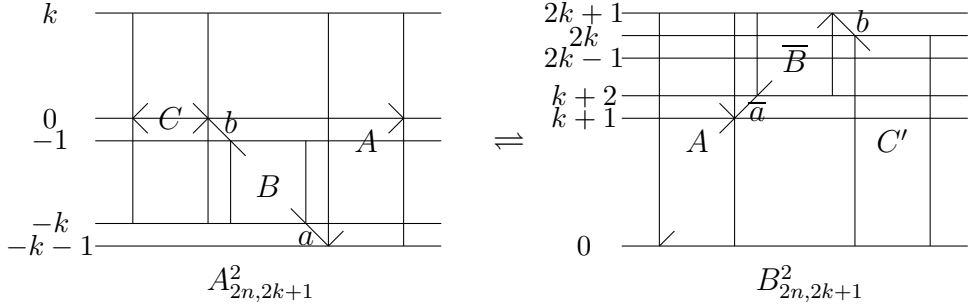

First, we have the bijection $\alpha_{2n,2k}$ in Theorem \ref{2n-2k} between $A_{2n,2k+1}^1$ and $B_{2n,2k+1}^1$.
Next, we establish the bijection between $A_{2n,2k+1}^2$ and $B_{2n,2k+1}^2$ in Figure \ref{2n-2k+1-2}. Since the process of the bijection is similar to that in Theorem \ref{2n-2k}, we only explain some special steps here. For the path on the left in Figure \ref{2n-2k+1-2}, the step $a$ is the first step ending on height $-k-1$. Then from $a$, we go back from right to left until we first arrive at height $0$. Let $b$ denote this step. For the path on the right in Figure \ref{2n-2k+1-2}, the step $b$ is the last down step ending on height $2k$. Then from $b$, we go back from right to left until we first arrive at height $k+1$, and let $\overline{a}$ denote this step. Notice that to obtain $C'$ which corresponds to
$C$, we first apply the bijection in Theorem \ref{2n-2k} on $C$, and then flip the corresponding path upside down.

Therefore, combining the two bijections $\alpha_{2n,2k}$ and $\alpha_{2n,2k+1}^2$, we prove the theorem.
\qed

\begin{thm}\label{2n+1-2k}
There is a bijection between $A_{2n+1,2k}$ and $B_{2n+1,2k}$.
\end{thm}

\pf We prove the theorem by induction on $k$.

When $k=1$, we separate $A_{2n+1,2}$ (resp.\  $B_{2n+1,2}$) into two subsets $A_{2n+1,2}^1$ and $A_{2n+1,2}^2$ (resp.\  $B_{2n+1,2}^1$ and $B_{2n+1,2}^2$).

$A_{2n+1,2}^1$ (resp.\  $A_{2n+1,2}^2$): the set of  paths in $A_{2n+1,2}$ which start with an up (resp.\  a down) step;

$B_{2n+1,2}^1$ (resp.\  $B_{2n+1,2}^2$): the set of  paths in $B_{2n+1,2}$ which end with a down (resp.\  an up ) step.

\begin{figure}[h]
\begin{center}
\begin{tikzpicture}


\draw(0,5)--(4,5);
\draw(0,4)--(4,4);
\draw(0,3)--(4,3);
\node at (-0.3,5){$1$};
\node at (-0.3,4){$0$};
\node at (-0.3,3){$-1$};

\draw(1,4)--(2,5);
\node at (1.3,4.6){$a$};

\draw(2,3)--(2,5);
\draw(3,3)--(3,5);
\draw(2,5)--(2.2,4.8);
\draw(3,3)--(2.8,3.2);
\node at (2.5,4){$A$};

\node at (2,2.5){$A_{2n+1,2}^1$};

\node at (4.8,4){$\rightleftharpoons$};

\draw(6,5)--(10,5);
\draw(6,4)--(10,4);
\draw(6,3)--(10,3);
\node at (5.8,5){$2$};
\node at (5.8,4){$1$};
\node at (5.8,3){$0$};

\draw(7,3)--(7,5);
\draw(8,3)--(8,5);
\draw(7,3)--(7.2,3.2);
\draw(8,5)--(7.8,4.8);
\node at (7.5,4){$\overline{A}$};

\draw(8,5)--(9,4);
\node at (8.8,4.6){$\overline{a}$};

\node at (7.5,2.5){$B_{2n+1,2}^1$};


\draw(0,2)--(4,2);
\draw(0,1)--(4,1);
\draw(0,0)--(4,0);
\node at (-0.3,2){$1$};
\node at (-0.3,1){$0$};
\node at (-0.3,0){$-1$};

\draw(1,1)--(2,0);
\node at (1.3,0.4){$a$};

\draw(2,0)--(2,2);
\draw(3,0)--(3,2);
\draw(2,0)--(2.2,0.2);
\draw(3,0)--(2.8,0.2);
\node at (2.5,1){$A$};

\node at (2,-0.5){$A_{2n+1,2}^2$};

\node at (4.8,1){$\rightleftharpoons$};

\draw(6,2)--(10,2);
\draw(6,1)--(10,1);
\draw(6,0)--(10,0);
\node at (5.8,2){$2$};
\node at (5.8,1){$1$};
\node at (5.8,0){$0$};

\draw(7,0)--(7,2);
\draw(8,0)--(8,2);
\draw(7,0)--(7.2,0.2);
\draw(8,0)--(7.8,0.2);
\node at (7.5,1){$A$};

\draw(8,0)--(9,1);
\node at (8.8,0.5){$\overline{a}$};

\node at (7.5,-0.5){$B_{2n+1,2}^2$};

\end{tikzpicture}
\end{center}
\caption{The bijection $\alpha_{2n+1,2}$ between $A_{2n+1,2}$ and $B_{2n+1,2}$}\label{2n+1-2}
\end{figure}
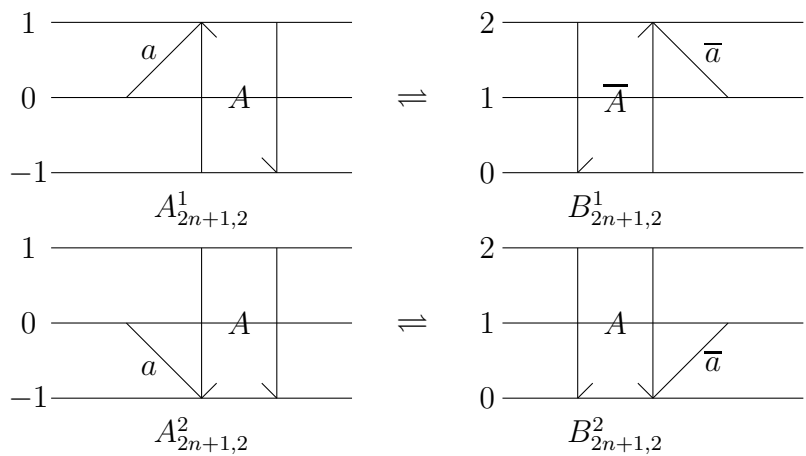

The bijections between $A_{2n+1,2}^i$ and $B_{2n+1,2}^i$ for $i=1,2$ are presented in Figure \ref{2n+1-2}.

Suppose that there is a bijection $\alpha_{2n+1,2k-2}$ between $A_{2n+1,2k-2}$ and $B_{2n+1,2k-2}$ ($k \geq 2$).
Let $A_{2n+1,2k}'$ denote the set of  paths in $A_{2n+1,2k}$ which touch the lines $y=k$ or $y=-k$.
Let $B_{2n+1,2k}'$ denote the set of  paths in $B_{2n+1,2k}$ such that in each path there exist some steps above the line  $y=2k-2$.
Since $A_{2n+1,2k}$ (resp.\  $B_{2n+1,2k}$) can be divided into $A_{2n+1,2k-2}$ and $A_{2n+1,2k}'$ (resp.\  $B_{2n+1,2k-2}$ and $B_{2n+1,2k}'$), we only need to build the bijection between $A_{2n+1,2k}'$ and $B_{2n+1,2k}'$. First, we separate $A_{2n+1,2k}'$ (resp.\  $B_{2n+1,2k}'$) into four subsets $A_{2n+1,2k}^1$, $A_{2n+1,2k}^2$, $A_{2n+1,2k}^3$, and $A_{2n+1,2k}^4$
(resp.\   $B_{2n+1,2k}^1$, $B_{2n+1,2k}^2$, $B_{2n+1,2k}^3$, and $B_{2n+1,2k}^4$).

$A_{2n+1,2k}^1$ (resp.\  $A_{2n+1,2k}^2$): the set of  paths in $A_{2n+1,2k}'$ which first arrive at height $k$, not height $-k$, such that all the steps (resp.\  some steps) before the first point on height $k$ are above (resp.\  below) height $0$.

$A_{2n+1,2k}^3$ (resp.\  $A_{2n+1,2k}^4$): the set of  paths in $A_{2n+1,2k}'$ which first arrive at height $-k$, not height $k$, such that all the steps (resp.\  some steps) before the first point on height $-k$ are below (resp.\  above) height $0$.

$B_{2n+1,2k}^1$ (resp.\  $B_{2n+1,2k}^3$): the set of  paths in $B_{2n+1,2k}'$ which end on height $2k-1$ such that the sub-path between the last up step starting on height $k-1$ and the last step arrives (resp.\  doesn't arrive) at height $2k$.

$B_{2n+1,2k}^2$ (resp.\  $B_{2n+1,2k}^4$): the set of  paths in $B_{2n+1,2k}'$ which end on heights $1,3, \ldots$, or $2k-3$. Especially, in each path, starting from the last down step ending on height $2k-2$, we go back from right to left until we arrive at height $k-1$ for the first time. This sub-path arrives (resp.\  doesn't arrive ) at height $2k$.

We state the bijections $\alpha_{2n+1,2k}^i$ between $A_{2n+1,2k}^i$ and $B_{2n+1,2k}^i$ for $i=1,2,3,4$ in Figures \ref{2n+1-2k-1}--\ref{2n+1-2k-4}, respectively. Since the bijections are similar to those in Theorem \ref{2n-2k}, we only explain some special steps. Notice that we still establish the bijections $\alpha_{2n+1,2k}^i$ for $i=2,4$ by induction on $k$.

For the path on the left in Figure \ref{2n+1-2k-1}, let $a$ denote the first up step ending on height $k$, and let $b$ denote the last down step ending on height $0$. On the other hand, for the path on the right in Figure \ref{2n+1-2k-1}, the last up step starting on height $k-1$ is denoted by $a$, and $\overline{b}$ denotes the the first up step ending on height $2k$ in the sub-path after the step $a$.

\begin{figure}[h]
\begin{center}
\begin{tikzpicture}


\draw(0,0)--(4.6,0);
\draw(0,0.8)--(4.6,0.8);
\draw(0,1.1)--(4.6,1.1);
\draw(0,1.4)--(4.6,1.4);
\draw(0,1.9)--(4.6,1.9);
\draw(0,2.2)--(4.6,2.2);
\node at (-0.6,0){\small$-k$};
\node at (-0.6,0.8){\small$-1$};
\node at (-0.6,1.1){\small$0$};
\node at (-0.6,1.4){\small$1$};
\node at (-0.6,1.9){\small$k-1$};
\node at (-0.6,2.2){\small$k$};

\draw(0.5,1.1)--(0.5,1.9);
\draw(1.5,1.1)--(1.5,1.9);
\draw(0.5,1.1)--(0.7,1.3);
\draw(1.5,1.9)--(1.3,1.7);
\node at (1,1.5){$B$};

\draw(1.5,1.9)--(1.8,2.2);
\node at (1.4,2.05){$a$};

\draw(1.8,2.2)--(1.8,0);
\draw(2.8,2.2)--(2.8,0);
\draw(1.8,2.2)--(2,1.9);
\draw(2.8,1.4)--(2.6,1.2);
\draw(2.8,1.4)--(2.6,1.6);
\node at (2.3,1.1){$A$};

\draw(2.8,1.4)--(3.1,1.1);
\node at (3.1,1.35){$b$};

\draw(3.1,1.1)--(3.1,0);
\draw(4.1,1.1)--(4.1,0);
\draw(3.1,1.1)--(3.3,0.9);
\draw(4.1,0.8)--(3.9,0.6);
\draw(4.1,0.8)--(3.9,1);
\node at (3.6,0.5){$C$};

\node at (2.2,-0.5){$A_{2n+1,2k}^1$};

\node at (5.4,1.1){$\rightleftharpoons$};

\draw(7,0)--(11.6,0);
\draw(7,0.8)--(11.6,0.8);
\draw(7,1.1)--(11.6,1.1);
\draw(7,1.6)--(11.6,1.6);
\draw(7,1.9)--(11.6,1.9);
\draw(7,2.2)--(11.6,2.2);
\node at (6.5,0){\small$0$};
\node at (6.5,0.8){\small$k-1$};
\node at (6.5,1.1){\small$k$};
\node at (6.5,1.6){\small$2k-2$};
\node at (6.5,1.9){\small$2k-1$};
\node at (6.5,2.2){\small$2k$};

\draw(7.5,0)--(7.5,2.2);
\draw(8.5,0)--(8.5,2.2);
\draw(7.5,0)--(7.7,0.2);
\draw(8.5,0.8)--(8.3,1);
\draw(8.5,0.8)--(8.3,0.6);
\node at (8,1.1){$\overline{A}$};

\draw(8.5,0.8)--(8.8,1.1);
\node at (8.85,0.95){$a$};

\draw(8.8,1.1)--(8.8,1.9);
\draw(9.8,1.1)--(9.8,1.9);
\draw(8.8,1.1)--(9,1.3);
\draw(9.8,1.9)--(9.6,1.7);
\node at (9.3,1.5){$B$};

\draw(9.8,1.9)--(10.1,2.2);
\node at (9.8,2.15){$\overline{b}$};

\draw(10.1,2.2)--(10.1,1.1);
\draw(11.1,2.2)--(11.1,1.1);
\draw(10.1,2.2)--(10.3,2);
\draw(11.1,1.9)--(10.9,1.7);
\draw(11.1,1.9)--(10.9,2.1);
\node at (10.6,1.6){$C$};

\node at (9.5,-0.5){$B_{2n+1,2k}^1$};

\end{tikzpicture}
\end{center}
\caption{The bijection $\alpha_{2n+1,2k}^1$ between $A_{2n+1,2k}^1$ and $B_{2n+1,2k}^1$}\label{2n+1-2k-1}
\end{figure}

For the path on the left in Figure \ref{2n+1-2k-2}, let $a$ denote the first up step ending on height $k$, and let $b$ denote the last down step ending on height $0$. Then from the step $a$, we go back from right to left until we first arrive at height $-1$, and let $c$ denote this step. For the path on the right in Figure \ref{2n+1-2k-2}, the last down step ending on height $2k-2$ is denoted by $b$. Then from $b$, we go back from right to left until we first arrive at height $k-1$, and let $a$ denote this step. Furthermore, in the sub-path between $a$ and $b$, from left to right, let $c$ denote the first up step ending on height $2k$. To obtain $D'$ which corresponds to $D$, we first use the hypothetical bijection on $D$. Then flip the corresponding path upside down.

\begin{figure}[h]
\begin{center}
\begin{tikzpicture}


\draw(0,0)--(5,0);
\draw(0,0.3)--(5,0.3);
\draw(0,0.8)--(5,0.8);
\draw(0,1.1)--(5,1.1);
\draw(0,1.4)--(5,1.4);
\draw(0,1.9)--(5,1.9);
\draw(0,2.2)--(5,2.2);
\node at (-0.6,0){\small$-k$};
\node at (-0.6,0.3){\small$-k+1$};
\node at (-0.6,0.8){\small$-1$};
\node at (-0.6,1.1){\small$0$};
\node at (-0.6,1.4){\small$1$};
\node at (-0.6,1.9){\small$k-1$};
\node at (-0.6,2.2){\small$k$};

\draw(0.5,0.3)--(0.5,1.9);
\draw(1.3,0.3)--(1.3,1.9);
\draw(0.5,1.1)--(0.7,0.9);
\draw(0.5,1.1)--(0.7,1.3);
\draw(1.3,0.8)--(1.1,1);
\draw(1.3,0.8)--(1.1,0.6);
\node at (0.9,1.1){$D$};

\draw(1.3,0.8)--(1.6,1.1);
\node at (1.6,0.95){$c$};

\draw(1.6,1.1)--(1.6,1.9);
\draw(2.4,1.1)--(2.4,1.9);
\draw(1.6,1.1)--(1.8,1.3);
\draw(2.4,1.9)--(2.2,1.7);
\node at (2,1.5){$B$};

\draw(2.4,1.9)--(2.7,2.2);
\node at (2.4,2.1){$a$};

\draw(2.7,2.2)--(2.7,0);
\draw(3.5,2.2)--(3.5,0);
\draw(2.7,2.2)--(2.9,2);
\draw(3.5,1.4)--(3.3,1.2);
\draw(3.5,1.4)--(3.3,1.6);
\node at (3.1,1.1){$A$};

\draw(3.5,1.4)--(3.8,1.1);
\node at (3.8,1.35){$b$};

\draw(3.8,1.1)--(3.8,0);
\draw(4.6,1.1)--(4.6,0);
\draw(3.8,1.1)--(4,0.9);
\draw(4.6,0.8)--(4.4,0.6);
\draw(4.6,0.8)--(4.4,1);
\node at (4.2,0.6){$C$};

\node at (2.5,-0.5){$A_{2n+1,2k}^2$};

\node at (5.5,1.1){$\rightleftharpoons$};

\draw(7,0)--(12,0);
\draw(7,0.8)--(12,0.8);
\draw(7,1.1)--(12,1.1);
\draw(7,1.6)--(12,1.6);
\draw(7,1.9)--(12,1.9);
\draw(7,2.2)--(12,2.2);
\node at (6.5,0){\small$0$};
\node at (6.5,0.8){\small$k-1$};
\node at (6.5,1.1){\small$k$};
\node at (6.5,1.6){\small$2k-2$};
\node at (6.5,1.9){\small$2k-1$};
\node at (6.5,2.2){\small$2k$};

\draw(7.5,0)--(7.5,2.2);
\draw(8.3,0)--(8.3,2.2);
\draw(7.5,0)--(7.7,0.2);
\draw(8.3,0.8)--(8.1,1);
\draw(8.3,0.8)--(8.1,0.6);
\node at (7.9,1.1){$\overline{A}$};

\draw(8.3,0.8)--(8.6,1.1);
\node at (8.6,0.95){$a$};

\draw(8.6,1.1)--(8.6,1.9);
\draw(9.4,1.1)--(9.4,1.9);
\draw(8.6,1.1)--(8.8,1.3);
\draw(9.4,1.9)--(9.2,1.7);
\node at (9,1.5){$B$};

\draw(9.4,1.9)--(9.7,2.2);
\node at (9.4,2.1){$c$};

\draw(9.7,2.2)--(9.7,1.1);
\draw(10.5,2.2)--(10.5,1.1);
\draw(9.7,2.2)--(9.9,2);
\draw(10.5,1.9)--(10.3,2.1);
\draw(10.5,1.9)--(10.3,1.7);
\node at (10.1,1.6){$C$};

\draw(10.5,1.9)--(10.8,1.6);
\node at (10.9,1.83){$b$};

\draw(10.8,1.6)--(10.8,0);
\draw(11.6,1.6)--(11.6,0);
\draw(10.8,1.6)--(11,1.4);
\node at (11.2,1){$D'$};

\node at (9.5,-0.5){$B_{2n+1,2k}^2$};

\end{tikzpicture}
\end{center}
\caption{The bijection $\alpha_{2n+1,2k}^2$ between $A_{2n+1,2k}^2$ and $B_{2n+1,2k}^2$}\label{2n+1-2k-2}
\end{figure}

For the path on the left in Figure \ref{2n+1-2k-3}, let $a$ denote the first down step ending on height $-k$. Alternatively, for the path on the right, the first up step starting on $k-1$ is denoted by $\overline{a}$.

\begin{figure}[h]
\begin{center}
\begin{tikzpicture}


\draw(0,0)--(3.3,0);
\draw(0,0.3)--(3.3,0.3);
\draw(0,0.8)--(3.3,0.8);
\draw(0,1.1)--(3.3,1.1);
\draw(0,2.2)--(3.3,2.2);
\node at (-0.6,0){\small$-k$};
\node at (-0.6,0.3){\small$-k+1$};
\node at (-0.6,0.8){\small$-1$};
\node at (-0.6,1.1){\small$0$};
\node at (-0.6,2.2){\small$k$};

\draw(0.5,0.3)--(0.5,1.1);
\draw(1.5,0.3)--(1.5,1.1);
\draw(0.5,1.1)--(0.7,0.9);
\draw(1.5,0.3)--(1.3,0.5);
\node at (1,0.7){$B$};

\draw(1.5,0.3)--(1.8,0);
\node at (1.47,0.13){$a$};

\draw(1.8,2.2)--(1.8,0);
\draw(2.8,2.2)--(2.8,0);
\draw(1.8,0)--(2,0.2);
\draw(2.8,0.8)--(2.6,1);
\draw(2.8,0.8)--(2.6,0.6);
\node at (2.3,1.1){$A$};

\node at (1.6,-0.5){$A_{2n+1,2k}^3$};

\node at (4,1.1){$\rightleftharpoons$};

\draw(5.5,0)--(8.8,0);
\draw(5.5,0.8)--(8.8,0.8);
\draw(5.5,1.1)--(8.8,1.1);
\draw(5.5,1.9)--(8.8,1.9);
\draw(5.5,2.2)--(8.8,2.2);
\node at (5,0){\small$0$};
\node at (5,0.8){\small$k-1$};
\node at (5,1.1){\small$k$};
\node at (5,1.9){\small$2k-1$};
\node at (5,2.2){\small$2k$};

\draw(6,0)--(6,2.2);
\draw(7,0)--(7,2.2);
\draw(6,0)--(6.2,0.2);
\draw(7,0.8)--(6.8,1);
\draw(7,0.8)--(6.8,0.6);
\node at (6.5,1.1){$A$};

\draw(7,0.8)--(7.3,1.1);
\node at (7.4,0.95){$\overline{a}$};

\draw(7.3,1.1)--(7.3,1.9);
\draw(8.3,1.1)--(8.3,1.9);
\draw(7.3,1.1)--(7.5,1.3);
\draw(8.3,1.9)--(8.1,1.7);
\node at (7.8,1.5){$\overline{B}$};

\node at (7,-0.5){$B_{2n+1,2k}^3$};

\end{tikzpicture}
\end{center}
\caption{The bijection $\alpha_{2n+1,2k}^3$ between $A_{2n+1,2k}^3$ and $B_{2n+1,2k}^3$}\label{2n+1-2k-3}
\end{figure}

For the path on the left in Figure \ref{2n+1-2k-4}, let $a$ denote the first down step ending on height $-k$, and let $b$ denote the last down step ending on height $0$ in the sub-path before the step $a$. For the path on the right in Figure \ref{2n+1-2k-4}, the last down step ending on height $2k-2$ is denoted by $b$. Then from $b$, we go back from right to left until we first arrive at height $k-1$, and let $\overline{a}$ denote this step. Furthermore, to obtain $C'$ which corresponds to $C$, we first use the hypothetical bijection on $\overline{C}$, and then flip the corresponding paths upside down.

\begin{figure}[h]
\begin{center}
\begin{tikzpicture}


\draw(0,0)--(4.6,0);
\draw(0,0.3)--(4.6,0.3);
\draw(0,0.8)--(4.6,0.8);
\draw(0,1.1)--(4.6,1.1);
\draw(0,1.4)--(4.6,1.4);
\draw(0,1.9)--(4.6,1.9);
\draw(0,2.2)--(4.6,2.2);
\node at (-0.6,0){\small$-k$};
\node at (-0.6,0.3){\small$-k+1$};
\node at (-0.6,0.8){\small$-1$};
\node at (-0.6,1.1){\small$0$};
\node at (-0.6,1.4){\small$1$};
\node at (-0.6,1.9){\small$k-1$};
\node at (-0.6,2.2){\small$k$};

\draw(0.5,0.3)--(0.5,1.9);
\draw(1.5,0.3)--(1.5,1.9);
\draw(0.5,1.1)--(0.7,1.3);
\draw(0.5,1.1)--(0.7,0.9);
\draw(1.5,1.4)--(1.3,1.2);
\draw(1.5,1.4)--(1.3,1.6);
\node at (1,1.1){$C$};

\draw(1.5,1.4)--(1.8,1.1);
\node at (1.8,1.35){$b$};

\draw(1.8,1.1)--(1.8,0.3);
\draw(2.8,1.1)--(2.8,0.3);
\draw(1.8,1.1)--(2,0.9);
\draw(2.8,0.3)--(2.6,0.5);
\node at (2.3,0.6){$B$};

\draw(2.8,0.3)--(3.1,0);
\node at (2.8,0.1){$a$};

\draw(3.1,2.2)--(3.1,0);
\draw(4.1,2.2)--(4.1,0);
\draw(3.1,0)--(3.3,0.2);
\draw(4.1,0.8)--(3.9,1);
\draw(4.1,0.8)--(3.9,0.6);
\node at (3.6,1.1){$A$};

\node at (2.2,-0.5){$A_{2n+1,2k}^4$};

\node at (5.4,1.1){$\rightleftharpoons$};

\draw(7,0)--(11.6,0);
\draw(7,0.8)--(11.6,0.8);
\draw(7,1.1)--(11.6,1.1);
\draw(7,1.6)--(11.6,1.6);
\draw(7,1.9)--(11.6,1.9);
\draw(7,2.2)--(11.6,2.2);
\node at (6.5,0){\small$0$};
\node at (6.5,0.8){\small$k-1$};
\node at (6.5,1.1){\small$k$};
\node at (6.5,1.6){\small$2k-2$};
\node at (6.5,1.9){\small$2k-1$};
\node at (6.5,2.2){\small$2k$};

\draw(7.5,0)--(7.5,2.2);
\draw(8.5,0)--(8.5,2.2);
\draw(7.5,0)--(7.7,0.2);
\draw(8.5,0.8)--(8.3,1);
\draw(8.5,0.8)--(8.3,0.6);
\node at (8,1.3){$A$};

\draw(8.5,0.8)--(8.8,1.1);
\node at (8.85,0.9){$\overline{a}$};

\draw(8.8,1.1)--(8.8,1.9);
\draw(9.8,1.1)--(9.8,1.9);
\draw(8.8,1.1)--(9,1.3);
\draw(9.8,1.9)--(9.6,1.7);
\node at (9.3,1.5){$\overline{B}$};

\draw(9.8,1.9)--(10.1,1.6);
\node at (10.15,1.85){$b$};

\draw(10.1,1.6)--(10.1,0);
\draw(11.1,1.6)--(11.1,0);
\draw(10.1,1.6)--(10.3,1.4);
\node at (10.6,0.8){$C'$};

\node at (9.5,-0.5){$B_{2n+1,2k}^4$};

\end{tikzpicture}
\end{center}
\caption{The bijection $\alpha_{2n+1,2k}^4$ between $A_{2n+1,2k}^4$ and $B_{2n+1,2k}^4$}\label{2n+1-2k-4}
\end{figure}
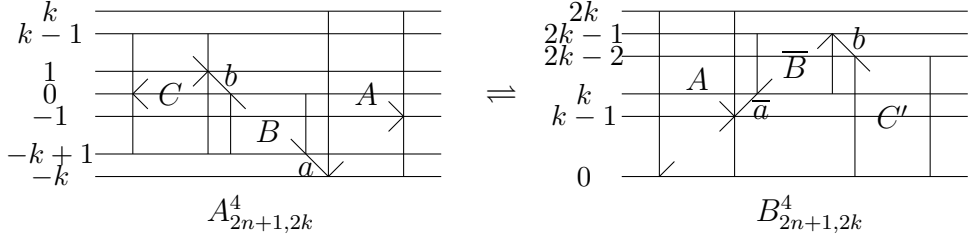

Notice that when $k=2$, the sub-path $D$ which belongs to the path on the left in Figure \ref{2n+1-2k-2} must end with a down step, and the sub-path $C$ belonging to the path on the left in Figure \ref{2n+1-2k-4} must end with an up step.

Combining the above bijections, we complete the proof.
\qed

\begin{thm}\label{2n+1-2k+1}
There is a bijection between $A_{2n+1,2k+1}$ and $B_{2n+1,2k+1}$.
\end{thm}
\pf  The set $A_{2n+1,2k+1}$ (resp.\  $B_{2n+1,2k+1}$) can be divided into three subsets $A_{2n+1,2k+1}^1$, $A_{2n+1,2k+1}^2$, and $A_{2n+1,2k+1}^3$
(resp.\  $B_{2n+1,2k+1}^1$, $B_{2n+1,2k+1}^2$, and $B_{2n+1,2k+1}^3$).

$A_{2n+1,2k+1}^1$: the set of  paths in $A_{2n+1,2k+1}$ which are in the strip $-k \leq y \leq k$.

$A_{2n+1,2k+1}^2$: the set of  paths in $A_{2n+1,2k+1}$ which arrive at height $-k-1$ such that all the steps before the first point on height $-k-1$ are below height $0$.

$A_{2n+1,2k+1}^3$: the set of  paths in $A_{2n+1,2k+1}$ which arrive at height $-k-1$ such that some steps before the first point on height $-k-1$ are above height $0$.

$B_{2n+1,2k+1}^1$: the set of  paths in $B_{2n+1,2k+1}$ which are in the strip $0 \leq y \leq 2k$.

$B_{2n+1,2k+1}^2$: the set of  paths in $B_{2n+1,2k+1}$ which end on height $2k+1$.

$B_{2n+1,2k+1}^3$: the set of  paths in $B_{2n+1,2k+1}$ which arrive at height $2k+1$, and end on heights $1,3,\ldots$, or $2k-1$.

\begin{figure}[h]
\begin{center}
\begin{tikzpicture}


\draw(0,0)--(3.2,0);
\draw(0,0.3)--(3.2,0.3);
\draw(0,0.8)--(3.2,0.8);
\draw(0,1.1)--(3.2,1.1);
\draw(0,1.9)--(3.2,1.9);
\node at (-0.6,0){\small$-k-1$};
\node at (-0.6,0.3){\small$-k$};
\node at (-0.6,0.8){\small$-1$};
\node at (-0.6,1.1){\small$0$};
\node at (-0.6,1.9){\small$k$};

\draw(0.5,0.3)--(0.5,1.1);
\draw(1.5,0.3)--(1.5,1.1);
\draw(0.5,1.1)--(0.7,0.9);
\draw(1.5,0.3)--(1.3,0.5);
\node at (1,0.7){$B$};

\draw(1.5,0.3)--(1.8,0);
\node at (1.45,0.15){$a$};

\draw(1.8,1.9)--(1.8,0);
\draw(2.8,1.9)--(2.8,0);
\draw(1.8,0)--(2,0.2);
\draw(2.8,0.8)--(2.6,0.6);
\draw(2.8,0.8)--(2.6,1);
\node at (2.3,0.9){$A$};

\node at (1.6,-0.5){$A_{2n+1,2k+1}^2$};

\node at (4.2,1){$\rightleftharpoons$};

\draw(6,0)--(9.3,0);
\draw(6,0.8)--(9.3,0.8);
\draw(6,1.1)--(9.3,1.1);
\draw(6,1.9)--(9.3,1.9);
\node at (5.5,0){\small$0$};
\node at (5.5,0.8){\small$k$};
\node at (5.5,1.1){\small$k+1$};
\node at (5.5,1.9){\small$2k+1$};

\draw(6.5,0)--(6.5,1.9);
\draw(7.5,0)--(7.5,1.9);
\draw(6.5,0)--(6.7,0.2);
\draw(7.5,0.8)--(7.3,0.6);
\draw(7.5,0.8)--(7.3,1);
\node at (7,0.9){$A$};

\draw(7.5,0.8)--(7.8,1.1);
\node at (7.85,0.9){$\overline{a}$};

\draw(7.8,1.1)--(7.8,1.9);
\draw(8.8,1.1)--(8.8,1.9);
\draw(7.8,1.1)--(8,1.3);
\draw(8.8,1.9)--(8.6,1.7);
\node at (8.3,1.5){$\overline{B}$};

\node at (7.6,-0.5){$B_{2n+1,2k+1}^2$};

\end{tikzpicture}
\end{center}
\caption{The bijection $\alpha_{2n+1,2k+1}^2$ between $A_{2n+1,2k+1}^2$ and $B_{2n+1,2k+1}^2$}\label{2n+1-2k+1-2}
\end{figure}
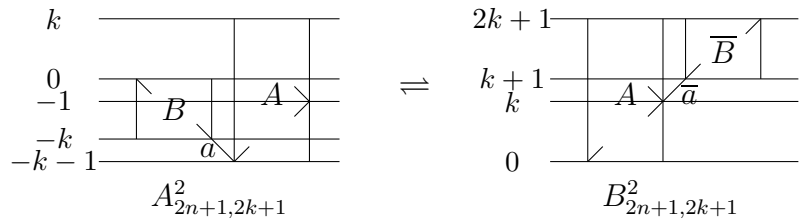

\begin{figure}[h]
\begin{center}
\begin{tikzpicture}


\draw(0,0)--(4.6,0);
\draw(0,0.3)--(4.6,0.3);
\draw(0,0.8)--(4.6,0.8);
\draw(0,1.1)--(4.6,1.1);
\draw(0,1.4)--(4.6,1.4);
\draw(0,1.9)--(4.6,1.9);
\node at (-0.6,0){\small$-k-1$};
\node at (-0.6,0.3){\small$-k$};
\node at (-0.6,0.8){\small$-1$};
\node at (-0.6,1.1){\small$0$};
\node at (-0.6,1.4){\small$1$};
\node at (-0.6,1.9){\small$k$};

\draw(0.5,0.3)--(0.5,1.9);
\draw(1.5,0.3)--(1.5,1.9);
\draw(0.5,1.1)--(0.7,0.9);
\draw(0.5,1.1)--(0.7,1.3);
\draw(1.5,1.4)--(1.3,1.2);
\draw(1.5,1.4)--(1.3,1.6);
\node at (1,0.95){$C$};

\draw(1.5,1.4)--(1.8,1.1);
\node at (1.85,1.3){$b$};

\draw(1.8,1.1)--(1.8,0.3);
\draw(2.8,1.1)--(2.8,0.3);
\draw(1.8,1.1)--(2,0.9);
\draw(2.8,0.3)--(2.6,0.5);
\node at (2.3,0.7){$B$};

\draw(2.8,0.3)--(3.1,0);
\node at (2.8,0.1){$a$};

\draw(3.1,1.9)--(3.1,0);
\draw(4.1,1.9)--(4.1,0);
\draw(3.1,0)--(3.3,0.2);
\draw(4.1,0.8)--(3.9,0.6);
\draw(4.1,0.8)--(3.9,1);
\node at (3.6,0.95){$A$};

\node at (2.2,-0.5){$A_{2n+1,2k+1}^3$};

\node at (5.4,1){$\rightleftharpoons$};

\draw(7,0)--(11.6,0);
\draw(7,0.8)--(11.6,0.8);
\draw(7,1.1)--(11.6,1.1);
\draw(7,1.6)--(11.6,1.6);
\draw(7,1.9)--(11.6,1.9);
\node at (6.5,0){\small$0$};
\node at (6.5,0.8){\small$k$};
\node at (6.5,1.1){\small$k+1$};
\node at (6.5,1.6){\small$2k$};
\node at (6.5,1.9){\small$2k+1$};

\draw(7.5,0)--(7.5,1.9);
\draw(8.5,0)--(8.5,1.9);
\draw(7.5,0)--(7.7,0.2);
\draw(8.5,0.8)--(8.3,0.6);
\draw(8.5,0.8)--(8.3,1);
\node at (8,0.95){$A$};

\draw(8.5,0.8)--(8.8,1.1);
\node at (8.85,0.95){$\overline{a}$};

\draw(8.8,1.1)--(8.8,1.9);
\draw(9.8,1.1)--(9.8,1.9);
\draw(8.8,1.1)--(9,1.3);
\draw(9.8,1.9)--(9.6,1.7);
\node at (9.3,1.5){$\overline{B}$};

\draw(9.8,1.9)--(10.1,1.6);
\node at (10.2,1.8){$b$};

\draw(10.1,1.6)--(10.1,0);
\draw(11.1,1.6)--(11.1,0);
\draw(10.1,1.6)--(10.3,1.4);
\node at (10.6,0.95){$C'$};

\node at (9.5,-0.5){$B_{2n+1,2k+1}^3$};

\end{tikzpicture}
\end{center}
\caption{The bijection $\alpha_{2n+1,2k+1}^3$ between $A_{2n+1,2k+1}^3$ and $B_{2n+1,2k+1}^3$}\label{2n+1-2k+1-3}
\end{figure}
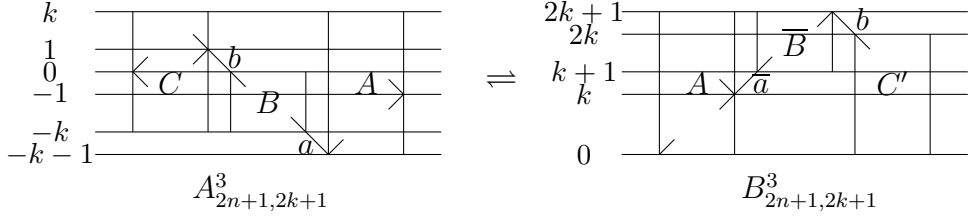

First, we have the bijection $\alpha_{2n+1,2k}$ in Theorem \ref{2n+1-2k} between $A_{2n+1,2k+1}^1$ and $B_{2n+1,2k+1}^1$.
Next, we build the bijections $\alpha_{2n+1,2k+1}^i$ between $A_{2n+1,2k+1}^i$ and $B_{2n+1,2k+1}^i$ for $i=2,3$ in Figures \ref{2n+1-2k+1-2} and \ref{2n+1-2k+1-3}, respectively. In the following, we again only explain some special steps.

For the path on the left in Figure \ref{2n+1-2k+1-2}, let $a$ denote the first down step ending on height $-k-1$. For the path on the right, the last up step starting on height $k$ is denoted by $\overline{a}$.

For the path on the left in Figure \ref{2n+1-2k+1-3}, let $a$ denote the first down step ending on height $-k-1$, and let $b$ denote the last down step ending on height $0$ in the sub-path before the step $a$. For the path on the right, the last down step ending on height $2k$ is denoted by $b$. Then from $b$, we go back from right to left until we first arrive at height $k$, and let $\overline{a}$ denote this step. Furthermore, to obtain $C'$ which corresponds to $C$, we first apply the bijection in Theorem \ref{2n+1-2k} on $\overline{C}$, and then flip the corresponding path upside down.

With the aid of the above bijections, we complete the proof.\qed

Therefore, combining Theorems \ref{2n-2k}--\ref{2n+1-2k+1}, we prove Theorem \ref{n-k}.

\section{Generating functions}
Although the emphasis of this paper is on bijections, we want to briefly demonstrate how the relevant identities can be seen via generating functions. We only do this for the instance $|A_{2n,2k}|=|B_{2n,2k}|$, but the other instances are similar.
 	
Lattice paths living in a strip have been treated in \cite{Panny-Prodinger}, so we can be brief. The series
 	\begin{equation*}
\varphi_0(z)=\sum_{n\ge0}|A_{n,2k}|z^n,
 	\end{equation*}
 	where $\varphi_0(z)$ is defined via the linear system (of $2k+1$ equations)
 	\begin{equation*}
\begin{bmatrix}
	1&-z&0\dots\\
-z&1&-z&0\dots\\
&&\dots\\
&&&	-z&1&-z\\
&&&&	-z&1\\
\end{bmatrix}
	\begin{bmatrix}
		\varphi_{-k}\\
		\varphi_{-k+1}\\
		\dots\\
		\varphi_{k-1}\\
\varphi_{k}\\
	\end{bmatrix}
	=\begin{bmatrix}
		0\\
				\dots\\
		1\\
		\dots\\
		0
	\end{bmatrix}.
 	\end{equation*}
 	The system can be solved using Cramer's rule, and one finds
 	\begin{equation*}
\varphi_{0}(z)=\frac{d_k^2}{d_{2k+1}},
 	\end{equation*}
 	where $d_j$ is the determinant of the system with $j$ rows and columns. Expanding, one finds the recursion $d_j=d_{j-1}-z^2d_{j-2}$, with the solution
 	\begin{equation*}
d_j=\frac{1-v^{2j+2}}{1-v^2}\frac1{(1+v^2)^j},
 	\end{equation*}
 	using the substitution $z=\frac{v}{1+v^2}$ for convenience.
 	
 	For the paths related to $B_{n,k}$, there is a similar system:
\begin{equation*}
	\begin{bmatrix}
		1&-z&0\dots\\
		-z&1&-z&0\dots\\
		&&\dots\\
		&&&	-z&1&-z\\
		&&&&	-z&1\\
	\end{bmatrix}
	\begin{bmatrix}
		\psi_{0}\\
		\psi_{1}\\
		\dots\\
		\psi_{2k-1}\\
		\psi_{2k}\\
	\end{bmatrix}
	=\begin{bmatrix}
		1\\
		\dots\\
		0\\
		\dots\\
		0
	\end{bmatrix}.
\end{equation*}
One finds with a similar argument that
 \begin{equation*}
\psi_j(z)=\frac{z^jd_{2k-j}}{d_{2k+1}},
 \end{equation*}	
so that we are left to show that
 \begin{equation*}
 	\frac{(1-v^{2k+2})^2}{(1-v^2)^2}\frac1{(1+v^2)^{2k}}=
d_k^2=\sum_{j=0}^kz^{2j}d_{2k-2j}=\sum_{j=0}^k\frac{v^{2j}}{(1+v^2)^{2j}}\frac{1-v^{4k-4j+2}}{1-v^2}
\frac1{(1+v^2)^{2k-2j}},
 \end{equation*}
or, simplified
 \begin{equation*}
 	\frac{(1-v^{2k+2})^2}{1-v^2}=\sum_{j=0}^k(v^{2j}-v^{4k-2j+2}),
 \end{equation*}
which is easy to check directly.

\section{Paths of width $3$ and bijections to certain families of trees}

We come back to the instance of $A_{n,3}$ and $B_{n,3}$ (which is already covered by our previous analysis),
since Cigler \cite{cigler} asked this question independently and since it leads to surprising bijections with other
mathematical objects.

We write $A(n,3|i)$ for the subfamily of $A_{n,3}$ ending on height $i$, and similarly for $B(n,3|i)$.

It is straightforward to prove that $|B(2n,3|0)|=F_{2n-1}$, $|B(2n,3|2)|=F_{2n}$, $|B(2n+1,3|1)|=F_{2n+1}$, and $|B(2n+1,3|3)|=F_{2n}$, with Fibonacci numbers $F_k$.

It is likewise straightforward to prove that $|A(2n,3|0)|=F_{2n+1}$, $|A(2n,3|-2)|=F_{2n}$, $|A(2n+1,3|1)|=F_{2n+1}$, and $|A(2n+1,3|=-1)|=F_{2n+2}$.

Thus we have that
\begin{equation*}
\Big|\bigcup_{i=0}^3B(n,3|i)\Big|=|A(n,3|0)|+|A(n,3|-1)|.
\end{equation*}

We are going  to link these problems to a tree structure named \emph{Elena trees},  introduced in \cite{elena}; compare also \cite{Deutsch-Prodinger}. The bijection presented in this early paper can be used to explain the equality.

If we want that the paths are in correspondence with trees, we require an even number of steps.

\subsection*{Paths in $B(2n,3|0)$ and height restricted plane trees}

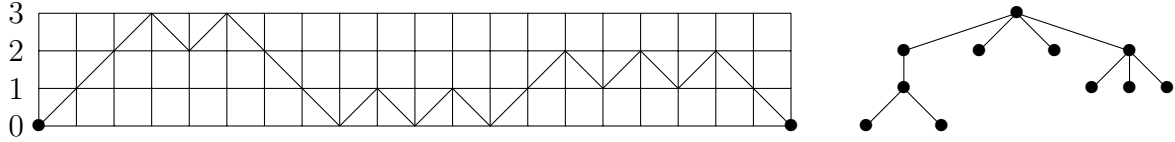
\begin{figure}[h]
	\begin{tikzpicture}[scale=0.5]
	\draw[step=1.cm,black] (0,0) grid (20,3);	
	\node at (-0.6,3) {$3$};	\node at (-0.6,2) {$2$};	\node at (-0.6,0) {$0$};\node at (-0.6,1) {$1$};
	\draw (0,0)to(3,3)to(4,2) to (5,3)to (8,0);
		\draw (8,0)to(9,1)to(10,0) to (11,1)to (12,0)to(14,2)to(15,1)to(16,2)to(17,1)to(18,2)to(20,0);

	\node at (0,0) {$\bullet$};\node at (20,0) {$\bullet$};
	
	\node at (26,4-1) {$\bullet$};
	\node at (23,3-1) {$\bullet$};
	\node at (25,3-1) {$\bullet$};
	\node at (27,3-1) {$\bullet$};
	\node at (29,3-1) {$\bullet$};
	\node at (23,2-1) {$\bullet$};
	\node at (29,2-1) {$\bullet$};
	\node at (28,2-1) {$\bullet$};
	\node at (30,2-1) {$\bullet$};
	\node at (22,1-1) {$\bullet$};
	\node at (24,1-1) {$\bullet$};

	\draw (22,0)to(23,1)to(23,2)to (26,3);\draw (24,0)to(23,1);\draw (25,2)to(26,3)to(27,2);
	\draw (26,3)to(29,2);
	\draw (28,1)to(29,2);
	\draw (29,1)to(29,2);
	\draw (30,1)to(29,2);
	
	\end{tikzpicture}
	\caption{A path of length 20, and the corresponding height restricted plane tree with 11 nodes}
\label{exam}
	\end{figure}

The translation of such a path of length $2n$ into a plane tree of height $\le 3$ (counting edges) is direct and sometimes called \emph{glove bijection}. The following example in Figure \ref{exam} will be sufficient.

\subsection*{Paths in $A(2n,3|0)$ and Elena trees}
	
	Elena trees were introduced in \cite{elena}; they consist of some nodes labelled $\mathbf{a}$, and a sequence of paths
	of various lengths (possibly empty)  emanating from all of them, except for the last one. An example in Figure \ref{Elena} describes this readily:
\begin{figure}[h]
\begin{tikzpicture}[scale=0.5]
		
\draw[ultra thick](0,0)to(4,-4);
		
		\node at (0,0) {$\bullet$};\node at (1,-1) {$\bullet$};\node at (2,-2) {$\bullet$};\node at (3,-3) {$\bullet$};
		\node at (4,-4) {$\bullet$};
		
		\draw(0,0)to(-4.5,-3);\node at (-1.5,-1) {$\bullet$};\node at (-3,-2) {$\bullet$};\node at (-4.5,-3) {$\bullet$};
		\draw(1,-1)to(-2,-2);\draw(1,-1)to(-1.0,-2);\draw(1,-1)to(-3,-5);
		\node at (-2,-2) {$\bullet$};\node at (-1,-2) {$\bullet$};\node at (0,-2) {$\bullet$};\node at (-1,-3) {$\bullet$};
		\node at (-2,-4) {$\bullet$};\node at (-3,-5) {$\bullet$};
		
		\draw(3,-3)to(1,-5);\node at (1,-5) {$\bullet$};\node at (2,-4) {$\bullet$};
		
		\node at (4.5,-3.5){$\mathbf{a}$};
				\node at (3.5,-2.5){$\mathbf{a}$};
						\node at (2.5,-1.5){$\mathbf{a}$};
								\node at (1.5,-0.5){$\mathbf{a}$};
										\node at (0.5,0.5){$\mathbf{a}$};
		\end{tikzpicture}
		\caption{An Elena tree described by  $\mathbf{a}\mathbf{p}_3\mathbf{a}\mathbf{p}_1\mathbf{p}_1\mathbf{p}_4\mathbf{a}\mathbf{a}\mathbf{p}_2	\mathbf{a}$}
\label{Elena}
\end{figure}
		
Typically, an Elena tree can be described by $\mathbf{a}\mathbf{p}_{i_1}\mathbf{p}_{i_2}\dots \mathbf{a}\mathbf{p}_{j_1}\mathbf{p}_{j_2}\dots \mathbf{a} \dots \mathbf{a}$.
For the set (language) of Elena trees, we might write a symbolic expression $(\mathbf{a}\mathbf{p}^*)^*\mathbf{a}$.
		
It is perhaps surprising that the paths in $A(2n,3|0)$ are suitable to describe Elena trees.
		For each sequence of steps $(2i,0)\to(2i+1,1)\to(2i+2,0)$, we write a symbol $\mathbf{a}$. In Figure \ref{p28}
		such pairs of steps are depicted in boldface.
		
		Thus, a path can be decomposed as $\mathbf{w}_0\mathbf{a}\mathbf{w}_1\mathbf{a}\dots \mathbf{a}\mathbf{w}_s$, where each $\mathbf{w}$ is a path from level 0 to level $0$ that
	``lives'' on levels $0,-1,-2$. Now we add a symbol $\mathbf{a}$ both, to the left and to the right.
	
	What is still left to be seen is how such a $\mathbf{w}$ can be interpreted as a sequence of paths: Each return to the level $0$ marks the end of a path, and the translation of the sojourns is as follows:
	
	$\tikz [scale=0.2]\draw[thick] (0,0)to (1,-1)to (2,0);$ corresponds to $\mathbf{p}_1$,
	$\tikz [scale=0.2]\draw[thick] (0,0)to (2,-2)to (4,0);$ corresponds to $\mathbf{p}_2$,
	$\tikz [scale=0.2]\draw[thick] (0,0)to (2,-2)to (3,-1) to (4,-2)to (6,0);$ corresponds to $\mathbf{p}_3$,
		$\tikz [scale=0.2]\draw[thick] (0,0)to (2,-2)to (3,-1) to (4,-2)to (5,-1) to (6,-2)to (8,0);$ corresponds to $\mathbf{p}_4$, and so forth.

Note that in this way a path of length $2n$ is (bijectively) mapped to an Elena tree of size (= number of nodes) $n+2$; the Elena tree consisting only of one node will not be considered.
		
\begin{figure}[h]
	\begin{tikzpicture}[scale=0.5]
	\draw[step=1.cm,black] (0,1) grid (28,-2);	
		\node at (-0.6,-1) {$-1$};	\node at (-0.6,-2) {$-2$};	\node at (-0.6,0) {$0$};\node at (-0.6,1) {$1$};
			\draw (0,0)to(2,-2)to(3,-1) to (4,-2)to (6,0);
			\draw[ultra thick](6,0)to(7,1)to(8,0);
			\draw (8,0)to(9,-1)to(10,0) to (11,-1)to (12,0)to(14,-2)to(15,-1)to(16,-2)to(17,-1)to(18,-2)to(20,0);
			\draw[ultra thick](20,0)to(21,1)to(22,0)to(23,1)to(24,0);
			\draw (24,0)to(26,-2)to(28,0);
			
			\node at (0,0) {$\bullet$};\node at (28,0) {$\bullet$};
		\end{tikzpicture}
			\caption{A path of length 28, described by $\mathbf{p}_3\mathbf{a}\mathbf{p}_1\mathbf{p}_1\mathbf{p}_4\mathbf{a}\mathbf{a}\mathbf{p}_2$}
\label{p28}
\end{figure}

\begin{figure}[h]
	\begin{tikzpicture}[scale=0.5]
	
	\draw[ultra thick](0,0)to(4,-4);
	
	\node at (0,0) {$\bullet$};\node at (1,-1) {$\bullet$};\node at (2,-2) {$\bullet$};\node at (3,-3) {$\bullet$};
	\node at (4,-4) {$\bullet$};
	
	\draw(0,0)to(-4.5,-3);\node at (-1.5,-1) {$\bullet$};\node at (-3,-2) {$\bullet$};\node at (-4.5,-3) {$\bullet$};
	\draw(1,-1)to(-2,-2);\draw(1,-1)to(-1.0,-2);\draw(1,-1)to(-3,-5);
	\node at (-2,-2) {$\bullet$};\node at (-1,-2) {$\bullet$};\node at (0,-2) {$\bullet$};\node at (-1,-3) {$\bullet$};
	\node at (-2,-4) {$\bullet$};\node at (-3,-5) {$\bullet$};
	
	\draw(3,-3)to(1,-5);\node at (1,-5) {$\bullet$};\node at (2,-4) {$\bullet$};
	\end{tikzpicture}
	\caption{The Elena tree with 16 nodes corresponding to  $(\mathbf{a})\mathbf{p}_3\mathbf{a}\mathbf{p}_1\mathbf{p}_1\mathbf{p}_4\mathbf{a}\mathbf{a}\mathbf{p}_2(\mathbf{a})$}
\end{figure}

\subsection*{Elena trees and height restricted plane trees}

We will establish a bijection between $A(2n,3|0)$ and
$B(2n,3|0)\cup B(2n,3|2)$; note, however that the latter set may be replaced by
$B(2n+2,3|0)$, by distinguishing the two cases of the last two steps.

So we would be done once we would know how to map (bijectively)
an Elena tree of size $n+2$ to a height restricted plane tree of the same size.

 This was documented already in \cite{elena}, but will be repeated here to make this description  self contained. The set of operations will be described as a sequence of pictures, which require no additional explanation.

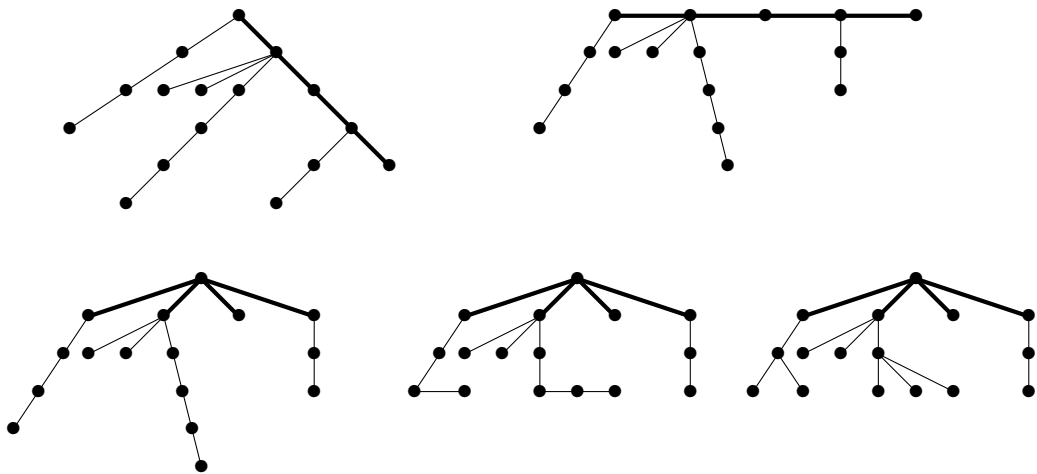
\begin{figure}[h]
 	\begin{tikzpicture}[scale=0.5]
 	
 	\draw[ultra thick](0,0)to(4,-4);
 	
 	\node at (0,0) {$\bullet$};\node at (1,-1) {$\bullet$};\node at (2,-2) {$\bullet$};\node at (3,-3) {$\bullet$};
 	\node at (4,-4) {$\bullet$};
 	
 	\draw(0,0)to(-4.5,-3);\node at (-1.5,-1) {$\bullet$};\node at (-3,-2) {$\bullet$};\node at (-4.5,-3) {$\bullet$};
 	\draw(1,-1)to(-2,-2);\draw(1,-1)to(-1.0,-2);\draw(1,-1)to(-3,-5);
 	\node at (-2,-2) {$\bullet$};\node at (-1,-2) {$\bullet$};\node at (0,-2) {$\bullet$};\node at (-1,-3) {$\bullet$};
 	\node at (-2,-4) {$\bullet$};\node at (-3,-5) {$\bullet$};
 	
 	\draw(3,-3)to(1,-5);\node at (1,-5) {$\bullet$};\node at (2,-4) {$\bullet$};

 \begin{scope}[shift={(10,0)}]
 	\draw[ultra thick](0,0)to(8,0);
 	
 	\node at (0,0) {$\bullet$};\node at (2,0) {$\bullet$};\node at (4,0) {$\bullet$};\node at (6,0) {$\bullet$};
 	\node at (8,0) {$\bullet$};
 	
 	\draw(0,0)to(-2,-3);\node at (-2/3,-1) {$\bullet$};\node at (-4/3,-2) {$\bullet$};\node at (-2,-3) {$\bullet$};
 	\draw(2,0)to(1,-1);\draw(2,0)to(0,-1);\draw(2,0)to(3,-4);
 	\node at (3,-4) {$\bullet$};
 	\node at (2.5,-2) {$\bullet$};\node at (2.25,-1) {$\bullet$};\node at (2.75,-3) {$\bullet$};
 	\node at (0,-1) {$\bullet$};\node at (1,-1) {$\bullet$};
 	
 	\draw(6,0)to(6,-2);\node at (6,-2) {$\bullet$};\node at (6,-1) {$\bullet$};

 \end{scope}
 \begin{scope}[shift={(-4,-8)}]
 \draw[ultra thick](0,0)to(3,1);
 \draw[ultra thick](2,0)to(3,1);
  \draw[ultra thick](4,0)to(3,1);
   \draw[ultra thick](6,0)to(3,1);

 \node at (0,0) {$\bullet$};\node at (2,0) {$\bullet$};\node at (4,0) {$\bullet$};\node at (6,0) {$\bullet$};
 \node at (3,1) {$\bullet$};

 \draw(0,0)to(-2,-3);\node at (-2/3,-1) {$\bullet$};\node at (-4/3,-2) {$\bullet$};\node at (-2,-3) {$\bullet$};
 \draw(2,0)to(1,-1);\draw(2,0)to(0,-1);\draw(2,0)to(3,-4);
 \node at (3,-4) {$\bullet$};
 \node at (2.5,-2) {$\bullet$};\node at (2.25,-1) {$\bullet$};\node at (2.75,-3) {$\bullet$};
 \node at (0,-1) {$\bullet$};\node at (1,-1) {$\bullet$};

 \draw(6,0)to(6,-2);\node at (6,-2) {$\bullet$};\node at (6,-1) {$\bullet$};

 \end{scope}
 \begin{scope}[shift={(6,-8)}]
 \draw[ultra thick](0,0)to(3,1);
 \draw[ultra thick](2,0)to(3,1);
 \draw[ultra thick](4,0)to(3,1);
 \draw[ultra thick](6,0)to(3,1);

 \node at (0,0) {$\bullet$};\node at (2,0) {$\bullet$};\node at (4,0) {$\bullet$};\node at (6,0) {$\bullet$};
 \node at (3,1) {$\bullet$};

 \draw(0,0)to(-4/3,-2)to (0,-2);\node at (-2/3,-1) {$\bullet$};\node at (-4/3,-2) {$\bullet$};

 \node at (0,-2) {$\bullet$};
 \draw(2,0)to(1,-1);\draw(2,0)to(0,-1);\draw(2,0)to(2,-2)to(4,-2);
 \node at (4,-2) {$\bullet$};
 \node at (2,-2) {$\bullet$};\node at (2,-1) {$\bullet$};\node at (3,-2) {$\bullet$};
 \node at (0,-1) {$\bullet$};\node at (1,-1) {$\bullet$};

 \draw(6,0)to(6,-2);\node at (6,-2) {$\bullet$};\node at (6,-1) {$\bullet$};

 \end{scope}
 \begin{scope}[shift={(15,-8)}]
 \draw[ultra thick](0,0)to(3,1);
 \draw[ultra thick](2,0)to(3,1);
 \draw[ultra thick](4,0)to(3,1);
 \draw[ultra thick](6,0)to(3,1);

 \node at (0,0) {$\bullet$};\node at (2,0) {$\bullet$};\node at (4,0) {$\bullet$};\node at (6,0) {$\bullet$};
 \node at (3,1) {$\bullet$};

 \draw(0,0)to(-4/3,-2);\draw(-2/3,-1)to(0,-2);

 \node at (-2/3,-1) {$\bullet$};\node at (-4/3,-2) {$\bullet$};

 \node at (0,-2) {$\bullet$};
 \draw(2,0)to(1,-1);\draw(2,0)to(0,-1);\draw(2,0)to(2,-2);
 \draw(2,-1)to(3,-2);
\draw(2,-1)to(4,-2);

 \node at (4,-2) {$\bullet$};
 \node at (2,-2) {$\bullet$};\node at (2,-1) {$\bullet$};\node at (3,-2) {$\bullet$};
 \node at (0,-1) {$\bullet$};\node at (1,-1) {$\bullet$};

 \draw(6,0)to(6,-2);\node at (6,-2) {$\bullet$};\node at (6,-1) {$\bullet$};

 \end{scope}

 	\end{tikzpicture}

 \caption{Transforming an Elena tree into a height restricted plane tree}
 \label{e16}

\end{figure}

We start with our running example of an Elena tree of size 16 and gradually transform it in Figure \ref{e16}.
 			
\subsection*{Paths with an odd number of steps}
 	
Let us consider $A(2n-1,3|-1)$, enumerated by $F_{2n}$. If we augment one up-step at the end, we have Elena trees, but with the special property that the last group of paths is non-empty.
 	
One the other hand, if we consider $B(2n-1,3|1)\cup B(2n-1,3|3)$, which is equivalent to $B(2n,3|2)$, then we augment it with 2 down-steps. The resulting height restricted tree has the property that the rightmost leaf is on a level $\ge2$.
 	
A short reflection convinces us that the bijection described earlier also works bijectively on the two respective subclasses.

\noindent {\bf Acknowledgements:} The first author was supported by the National Natural Science Foundation of China and the Fundamental Research Funds for the Central Universities.

\end{document}